\documentclass[a4paper]{article}

\usepackage{latexsym}

\newtheorem{thm}{Theorem}[section]
\newtheorem{pro}[thm]{Proposition}
\newtheorem{lem}[thm]{Lemma}
\newtheorem{cor}[thm]{Corollary}

\def\qed{{$\Box$}\medskip}

\def\implies{\Rightarrow}

 \def\a{\alpha}   \def\d{\delta}
\def\e{\varepsilon} \def\f{\varphi}  
\def\p{\psi} \def\r{\rho} \def\s{\sigma}  
\def\om{\omega}   \def\D{\Delta}

 \def\bC{{\mathbf C}}  
  \def\bT{{\mathbf T}} 

 \def\Ad{\mathrm{Ad}}

 \def\ip#1#2{{(\,#1\, , \, #2\,)}}
 \def\n#1{\Vert #1 \Vert}
\def\demo{\noindent \emph{Proof.\ }}

\def\AA{{(A,\D)}} \def\aqg{algebraic quantum group}
\def\ar{A_{\mathrm{r}}} \def\dar{\Delta_{\mathrm{r}}}
\def\pr{\pi_{\mathrm{r}}}

\def\har{\hat A_{\mathrm{r}}}
\def\hdar{\hat \Delta_{\mathrm{r}}}
\def\hdarop{\hat \Delta_{\mathrm{r, op}}}

\def\au{A_{\mathrm{u}}} \def\dau{\Delta_{\mathrm{u}}}
\def\pu{\pi_{\mathrm{u}}} \def\eu{\varepsilon_{\mathrm{u}}}

\def\LL{\Lambda} \def\HAA{(\hat A,\hat \D)}
 
\def\drop{ \hat \Delta_{\mathrm{r,op}}}

\def\qr{Q_{\mathrm{r}}}

\def\F{{\cal F}}
\def\hs{\hat{S}_{\mathrm{r, op}}}
\def\hS{\hat{S}_{\mathrm{op}}}

\def\A{{\cal A}}

\def\implies{\Rightarrow}

\begin{document}

\title{Amenability and Co-Amenability of Algebraic Quantum Groups II}

\author{E.~Bedos$^*$\ \   G.J.~Murphy  \ \ L.~Tuset$^*$}
\date{October 2001}
\maketitle

\begin{abstract}
We continue our study of the concepts of amenability and co-amenability for
algebraic quantum groups in the sense of A. Van Daele
 and our investigation  of their relationship  with nuclearity
and injectivity. One major tool for our
analysis is that every non-degenerate *-representation of the universal
C*-algebra
associated to an \aqg\ has a unitary generator which may be described in
a concrete way.

\vspace{5ex}

\hspace{-\parindent} \emph{Subj. Class.}: Ouantum groups, C*-algebras

\hspace{-\parindent} \emph{MSC 2000}:  Primary 46L05, 46L65.
Secondary 16W30, 22D25, 58B32.

\hspace{-\parindent} \emph{Keywords}: quantum group, amenability,
nuclearity, injectivity.

\end{abstract}

\vfill
\thanks{$^*$  Partially supported by the Norwegian Research Council}

\newpage

\section{Introduction}

The present paper is a continuation of our previous
paper~\cite{BMT2}, where we initiated a study of the
concepts of amenability and co-amenability for algebraic quantum groups
(see  also \cite{BMT}).
We gave there several  equivalent
formulations of co-amenability
and showed that co-amenability of an \aqg\
$\AA$ always implies amenability of its dual algebraic quantum group
$(\hat A,\hat \D).$ We also obtained some results concerning
the relationship between co-amenability of
$\AA$, injectivity of the von Neumann algebra $M$ associated to $\AA$ and
amenability of $(\hat A,\hat \D)$. The  algebra $M$ is the von Neumann
algebra generated by $\ar$, where $(\ar,\dar)$ denotes the analytic
extension
of $\AA$.

One may  construct  a unique universal C$^{*}$-algebraic quantum group
$(A_{\mathrm{u}},\Delta_{\mathrm{u}})$ associated to an \aqg\ $\AA$ (see
\cite{Ku2}).
We  show in Section 4 of the present paper that co-amenability of $\AA$ is
equivalent to the
fact that the canonical homomorphism from $A_{\mathrm{u}}$ onto $\ar$ is
injective (see \cite{BMT} for the compact case). This generalizes a
classical result in the case that $\AA$ is the \aqg\  associated to the
group algebra
of a discrete group. Further, we  establish the
following result, which is also well-known in the
group algebra case (see \cite{Pa}) :

\begin{thm} \label{main} Let $\AA$ be an \aqg .

    Consider the following statements:

    \vspace{1ex}

\hspace{10ex}    (1) $\AA$ is co-amenable

\hspace{10ex}    (2) $(\hat A,\hat \D)$ is amenable

\hspace{10ex}    (3) $A_{u}$ is nuclear

\hspace{10ex}    (4) $\ar$ is nuclear

\hspace{10ex}    (5) $M$ is injective

   \vspace{1ex}
 Then $(1) \implies (2) \implies (3)
\implies (4) \implies (5)$.

\vspace{2ex}

If $\AA$ is compact and  has a tracial Haar functional, then we
also have \newline \indent  $(5)\implies (1)$, that is,  all statements
above are
equivalent.
\end{thm}

The main new  part of this result is the fact
that (2) implies
(3). It is possible to deduce that (1) implies (3) from C.-K. Ng's
paper \cite{Ng1} on Hopf C*-algebras. His proof is related to the
one given by  E. Blanchard
(\cite{Bl}, see also\cite{BS}) in the setting of regular
multiplicative unitaries. Our proof is quite different and relies  on the
characterization of the nuclearity of a C$^*$-algebra B in terms of the
injectivity of $ B^{**}.$
The equivalence between (1), (2) and (5) in the compact tracial case
may be deduced from Ruan's main result (Theorem 4.5) in
\cite{Ruan}. We propose a proof that (5)
implies (1) in this case, which we believe is somewhat more accessible than
his.
The interesting question as to whether any (or all) of the statements
(2),(3),(4) or (5)
always implies (1) seems quite hard to answer. As a pendant to this
question, we
 show
that for a compact $\AA,$ injectivity of $M$ always implies
a kind  of ``perturbed'' co-amenability, involving the notion of quantum
dimension (of irreducible unitary corepresentations).

An important tool in our approach is the fact that any non-degenerate
\\ $*$-representation
of $\au$ on some Hilbert space has a unitary ``generator'' which may
be described in a concrete way. We present a self-contained proof of this
Kirchberg-type
result
in  Section 3. Similar results (using the universal corepresentation
and a certain $L^1$-algebra) have been previously obtained  by J.
Kustermans in \cite{Ku2,Ku3}.  The classical result of  E. Kirchberg
for Kac algebras may be found in \cite{ES}.

The reader should consult \cite{BMT2} for an introduction to this
subject, including a more extensive list of references to related
papers. Section 2 contains a review of most of the necessary
background material required for understanding the present paper.

\section{Preliminaries}
We recall in this section some definitions and results from
\cite{BMT,KD,VD,VD2}. We also prove some technical lemmas that we need
later on.

\vspace{2ex}
We begin with some
terminology that will be used throughout the paper.

\vspace{1ex}
Every algebra will be a (not necessarily unital) associative
algebra over the complex field~$\bC$. The identity map on a
set $V$ will be denoted by $\iota_V$, or simply by $\iota$, if
no ambiguity is involved.

\vspace{1ex}
If $V$ and $W$ are linear spaces, $V'$ denotes the linear
space of linear functionals on $V$ and  $V \otimes W$ denotes
the linear space tensor product of $V$ and~$W$. The {\em flip
map}~$\chi$ from ${V \otimes W}$ to ${W \otimes V}$ is the
linear map sending ${v\otimes w}$ onto ${w\otimes v}$, for all
$v\in V$ and $w\in W$. If $V$ and $W$ are Hilbert spaces,
${V\otimes W}$ denotes their Hilbert space tensor product; we
denote by $B(V)$ and $B_0(V)$ the C*-algebras of bounded
linear operators and compact operators on~$V$, respectively.
If $v\in V$ and $w\in W$, $\om_{v,w}$ denotes the weakly
continuous  bounded linear functional on $B(V)$ that maps $x$
onto $\ip{x(v)}w$. We set $\om_v=\om_{v,v}$. We will often
also use the notation $\om_v$ to denote a restriction to a
C*-subalgebra of $B(V)$ (the domain of $\om_v$ will be
determined by the context).

\vspace{1ex}
If $V$ and $W$ are algebras, ${V \otimes W}$ denotes their
algebra tensor product. We sometimes denote this algebraic tensor product
by $V \odot W$ if we feel there is some danger of confusion.
If $V$ and $W$ are C*-algebras, then
${V\otimes W}$ will denote their C*-tensor product with
respect to the minimal (spatial) C*-norm. If $V$ and $W$ are von
Neumann algebras, then
${V \bar{\otimes} W}$ will denote their von Neumann algebra tensor product.

\vspace{2ex}
For a  review of some results related to multiplier algebras, especially
multiplier
algebras of C*-algebras, and  to  slice maps, we refer to
\cite{BMT}. We will use repeatedly these results and also most of the
terminology introduced in
this paper. For the ease of the reader and to fix notation, we  recall here
some of the
basic definitions and properties of algebraic quantum groups.

\vspace{2ex}
Let $A$ be a non-degenerate $*$-algebra and let $\D$ be a
non-degenerate $*$-homomorphism from $A$ into $M(A\otimes A)$.
Suppose that the following conditions hold:

\vspace{1ex}
(1) $(\D\otimes \iota)\D=(\iota\otimes \D)\D$;

(2) The linear mappings defined by the assignments
 ${a\otimes b\mapsto \D(a)(b\otimes 1)}$ and
\newline \indent ${a\otimes b\mapsto \D(a)(1\otimes b)}$ are bijections from
$A\otimes A$ onto itself.

\vspace{1ex}
\noindent Then the pair $(A,\D)$ is called a {\em multiplier Hopf
$*$-algebra}.

\vspace{2ex}
In Condition~(1), we are regarding both maps as maps into
$M(A\otimes A\otimes A)$, so that their equality makes sense.
It follows from Condition~(2), by taking adjoints, that the
maps defined by the assignments
 ${a\otimes b\mapsto (b\otimes 1)\D(a)}$ and
${a\otimes b\mapsto (1\otimes b)\D(a)}$ are also bijections
from $A\otimes A$ onto itself.

\vspace{2ex}
Let $\AA$ be a multiplier Hopf $*$-algebra and let $\omega$ be
a linear functional on $A$ and $a$ an element in $A$. There is
a unique element  $(\omega \otimes \iota) \Delta(a)$ in $M(A)$
for which \[(\omega \otimes \iota) (\Delta(a))b = (\omega
\otimes \iota)(\Delta(a)(1 \otimes b))\] and \[ b(\omega
\otimes \iota) (\Delta(a)) = (\omega \otimes \iota)((1 \otimes
b) \Delta(a)),\] for all $b \in A$. The element
$(\iota \otimes \omega) \Delta(a)$ in $M(A)$ is determined similarly.
Thus, $\om$ induces linear maps maps $(\om \otimes \iota)
\Delta$ and $(\iota \otimes \omega) \Delta$ from $A$ to
$M(A)$.

\vspace{2ex}
There exists a unique non-zero $*$-homomorphism $\varepsilon$
from $A$ to $\bf C$ such that, for all $a\in A$,
\[(\varepsilon \otimes \iota)\D(a) =(\iota \otimes
\varepsilon)\Delta(a) = a.\] The map $\varepsilon$ is called
the {\em co-unit} of $(A, \Delta)$. Also, there exists a
unique anti-multiplicative linear isomorphism $S$ on $A$ that
satisfies the conditions
\[m(S \otimes \iota) (\Delta(a)(1 \otimes b)) = \varepsilon
(a)b\] \noindent and
\[m(\iota \otimes S)((b\otimes 1)\Delta (a)) = \varepsilon (a)b,\]
\noindent for all $a,b \in A$. Here ${m\colon  A \otimes A\to
A}$ denotes the linearization  of the multiplication map
${A\times A\to A}$. The map $S$ is called the {\em antipode}
of $(A, \Delta)$.
The antipode is in general neither $*$-preserving, nor involutive;
however, we have $S(S(a^{*})^*) = a $ for all $a \in A.$

\vspace{2ex}

If $\om\in A'$, we say $\omega$ is {\em left
invariant} if $(\iota \otimes \omega ) \Delta(a) =
\omega(a)1$, for all $a \in A$. Right invariance is defined
similarly. If a non-zero left-invariant linear functional
on~$A$ exists, it is unique, up to multiplication by a
non-zero scalar. A similarl statement holds for a non-zero right-invariant
linear functional. If $\f$ is a left-invariant functional
on~$A$, the functional $\p=\f S$ is right invariant.

\vspace{2ex}
If $A$ admits a non-zero, left-invariant, positive linear
functional~$\f$, we call $(A, \Delta)$ an {\em algebraic
quantum group} and we call $\f$ a {\em left Haar integral}
on~$\AA$. Faithfulness of $\f$ is automatic.

\vspace{2ex}
Note that although $\p=\f S$ is right invariant, it may not be
positive. On the other hand, it is proved in \cite{KD} that a
non-zero, right-invariant, positive linear functional on
$A$---a {\em right Haar integral}---necessarily exists. As for
a left Haar integral, a right Haar integral is necessarily
faithful.

\vspace{2ex}

The left Haar functional $\f$ is not necessarily tracial (or central).
However, there is a unique bijective homomorphism $\r\colon A\to A$
such that $\f(a b)=\f(b \r(a))$, for all $a,b\in A$. Moreover,
$\r(\r(a^*)^*)=a$.

\vspace{2ex}

One useful property, which will be used several times in the sequel,
is that every element of $A$ has ``compact support'' : given $a \in
A$, there exists some $b \in A$ such that $a b = ba = a.$ ( In fact, a
more general result is true \cite{Ku1}.) Another property we will need
is the so-called strong left invariance of $\f,$ which is proved in
\cite[Proposition 3.11]{VD}. It says that
$$(\iota \odot \f)((1 \otimes a) \D(b)) = S((\iota \odot \f)( \D(a)(1
\otimes
b)))$$
holds for all $a, b \in A.$

\vspace{2ex} We now turn to a short discussion of duality.
If $\AA$ is an \aqg, denote by $\hat A$ the linear
subspace of $A'$ consisting  of all functionals $\varphi a$,
where  $a \in A$. Since $\f a=\r(a)\f$, we have $\hat A={\{a
\varphi\ | \ a \in A \}}$. If $\om_1,\om_2\in \hat A$, one can
define a linear functional ${(\om_1\otimes \om_2)\D}$ on~$A$
by setting ${(\om_1\otimes \om_2)\D(a)}={(\f\otimes
\f)((a_1\otimes a_2)\D(a))}$, where $\om_1=\f a_1$ and
$\om_2=\f a_2$. Using this,  the space $\hat A$ can be made
into a non-degenerate $*$-algebra. The multiplication is given
by $\omega_1 \omega_2 = (\omega_1 \otimes \omega_2)\Delta$ and
the involution is given by setting $\omega^*(a) =
\omega(S(a)^*)^-$, for all $a\in A$ and $\om_1,\om_2,\om\in
\hat A$; it is clear that $\omega_1 \omega_2, \om^*\in A'$ but
one can show that, in fact, $\om_1\om_2, \om^*\in \hat A$.

\vspace{2ex}
One can realize $M(\hat A)$ as a linear space by identifying
it as  the linear subspace of $A'$ consisting of all $\om\in
A'$ for which ${(\om\otimes \iota)\D(a)}$ and  ${(\iota\otimes
\om)\D(a)}$ belong to $A$. (It is clear that $\hat A$ belongs
to this subspace.) In this identification of $M(\hat A)$, the
multiplication and involution are determined by
\[ (\om_1\om_2)(a)=\om_1((\iota\otimes
\om_2)\D(a))=\om_2((\om_1\otimes \iota)\D(a))\] and
\[\om^*(a)=\om(S(a)^*)^-,\] for all $a\in A$ and
$\om_1,\om_2,\om\in M(\hat A)$.

Note that the co-unit $\e$ of $A$ is the unit of the
$*$-algebra~$M(\hat A)$.

\vspace{1ex}
There is a unique $*$-homomorphism  $\hat \Delta$ from $\hat
A$ to $M(\hat A\otimes \hat A)$ such that for all
$\om_1,\om_2\in \hat A$ and $a,b\in A$,
\[((\om_1\otimes 1)\hat \D(\om_2))(a\otimes b)=(\om_1\otimes
\om_2)(\D(a)(1\otimes b))\] and
\[(\hat \D(\om_1)(1 \otimes \om_2))(a\otimes b)=(\om_1\otimes
\om_2)((a\otimes 1)\D(b)).\] Of course, we are here
identifying ${A'\otimes A'}$ as a linear subspace of
${(A\otimes A)'}$ in the usual way, so that elements of ${\hat
A\otimes \hat A}$ can be regarded as linear functionals on
${A\otimes A}$.

\vspace{2ex}
The pair $(\hat A, \hat \Delta)$ is an algebraic quantum
group, called the {\em dual} of $\AA$. Its co-unit $\hat
\varepsilon$ and antipode $\hat S$ are given by $\hat
\varepsilon(a\varphi) = \varphi(a)$ and
$\hat S(a\varphi)= (a\varphi)\circ S$, for all $a\in A$.

\def\doubledual#1{#1\hat{}\phantom{,}\hat{}}
\vspace{2ex}

There is an algebraic quantum group version of Pontryagin's
duality theorem for locally compact abelian groups that
asserts that $\AA$ is canonically isomorphic to the dual of
$(\hat A, \hat \Delta)$; that is, $\AA$ is isomorphic to its
double dual $(\doubledual A,\doubledual \Delta)$.

\vspace{2ex}
We now turn to the analytic theory of algebraic quantum groups.
We first recall the concept of a GNS~pair. Suppose
given a positive linear functional $\omega$ on a
$*$-algebra~$A$. Let $H$ be a Hilbert space, and let
${\LL\colon A\to H}$ be a linear map with dense range for
which $\ip{\LL(a)}{\LL(b)}=\om(b^*a)$, for all $a,b\in A$.
Then we call $(H,\LL)$ a {\em GNS~pair} associated to~$\om$.
As is well known, such a pair always exists and is essentially
unique. For, if $(H',\LL')$ is another GNS~pair associated to
$\om$, the map, ${\LL(a)\mapsto \LL'(a)}$, extends to a
unitary $U:H\to H'$.

\vspace{2ex}
If $\f$ is a left Haar integral on an \aqg\ $\AA$, and
$(H,\LL)$ is an associated GNS~pair, then it can be shown that
there is a unique $*$-homomorphism ${\pi\colon A\to B(H)}$
such that $\pi(a)\LL(b)=\LL(a b)$, for all ${a,b\in A}$.
Moreover, $\pi$ is faithful and non-degenerate. We let $\ar $
denote the norm closure of $\pi(A)$ in $B(H)$. Thus, $\ar $ is
a non-degenerate C*-subalgebra of~$B(H)$. The $*$-representation
${\pi\colon A\to B(H)}$  is essentially unique, for if
$(H',\LL')$ is another GNS~pair associated to $\f$, and
${\pi'\colon A\to B(H')}$ is the corresponding $*$-representation,
then, as we observed above, there exists a unitary ${U\colon
H\to H'}$ such that $U\LL(a)=\LL'(a)$,  for all $a\in A$,  and
consequently, $\pi'(a)=U\pi(a)U^*$.

We shall use the symbol $M$ to denote the
von~Neumann algebra generated by~$\ar $. Of course, $\ar $ and
$\pi(A)$ are weakly dense in~$M$.

\vspace{2ex}
Now observe that there exists a unique non-degenerate
$*$-homomorphism ${\dar  : \ar  \rightarrow M(\ar  \otimes \ar
)}$ such that, for all $a \in A$ and  all $x \in A \otimes A$,
we have \[{\dar (\pi(a)) (\pi \otimes \pi) (x)} = {(\pi
\otimes \pi) (\Delta (a)x)}\] and \[{(\pi \otimes \pi) (x)\dar
(\pi(a))} = {(\pi \otimes \pi) (x \Delta (a))}.\]

We also recall that
\[ \ar = [(\om \otimes \iota)( \dar (x) ) \, | \ x \in \ar, \ \om \in
\ar^{*}] =  [(\iota \otimes \om)( \dar (x) ) \, | \ x \in \ar,
\ \om \in \ar^{*}]. \]

The pair $(\ar, \dar )$ is a reduced locally compact quantum
group in the sense of Definition 4.1 of \cite{KV};  it is called
the {\em analytic extension} of $\AA$ associated to~$\f$.

\vspace{2ex}
We also need to recall that there is a unique unitary operator
$W$ on $H \otimes H$ such that
\[W((\Lambda\otimes\Lambda)(\Delta(b)(a \otimes 1)))=
\Lambda (a) \otimes \Lambda (b),\] for all $a, b \in A$. This
unitary satisfies the equation \[W_{12} W_{13} W_{23} =
W_{23} W_{12};\] thus, it is a multiplicative unitary,
said to be {\em associated to}~$(H,\LL)$.
Here we have used the leg numbering notation of \cite{BS}.
One can show that
$W\in M(\ar \otimes B_0(H))$, so especially $W\in (\ar \otimes B_0(H))^{''}
= M \bar{\otimes} B(H)$, and that $\ar $ is the norm
closure of the linear space $\{(\iota \otimes \omega) (W)\mid
{\omega \in B_0(H)^*}\}$. Also,
 $\dar  (a) = {W^* (1 \otimes a) W}$,
for all  $a \in \ar $.

\vspace{2ex}

Since the map $\dar$ is
unitarily implemented, it has a unique weakly continuous
extension to a unital $*$-homomorphism ${\dar :M \rightarrow
M\bar{\otimes} M}$, given explicitly by $\dar (a )=W^*
(1\otimes a)W$, for all $a\in M $. The Banach space $M_{*}$ may
be regarded as a Banach
algebra when equipped  with the  canonical multiplication
induced by $\dar$; thus, the product of two elements $\om$ and
$\s$ is given by $\om \s={(\om \bar{\otimes}\s)\circ \dar}$.
 
We use the same symbol $R$ to denote the anti-unitary antipode
of $\ar$ and of~$M$, and we denote by $\tau$ the scaling group
of $(\ar,\dar)$ (see \cite{KD,KV}).

\vspace{2ex}
Consider now the algebraic dual $(\hat{A}, \hat{\D })$ of
$(A,\D )$. A right-invariant linear functional $\hat{\psi}$ is
defined on $\hat{A}$ by setting
$\hat{\psi}(\hat{a})=\varepsilon (a)$, for all $a\in A$. Here
$\hat{a}=a\varphi$ and $\e$ is the co-unit of~$\AA$. Since the
linear map, $A\rightarrow\hat{A}$, $a\mapsto\hat{a}$, is a
bijection (by faithfulness of  $\varphi$), the functional
$\hat{\psi}$ is well defined. Now define a linear map
$\hat{\Lambda} :\hat{A}\rightarrow H$ by setting
$\hat{\Lambda}(\hat {a})=\Lambda (a)$, for all $a\in\A$. Since
$\hat{\psi} (\hat{b}^*\hat{a})=\varphi (b^*
a)=\ip{\LL(a)}{\LL(b)}$, for all $a,b\in A$, it follows that
$(H,\hat{\Lambda} )$ is a GNS-pair associated to $\hat{\psi}$.
It can be shown that it is unitarily equivalent to the
GNS-pair for a  left Haar integral $\hat{\varphi}$ of $(\hat
A,\hat \D)$. Hence, we can use $(H,\hat\LL)$ to define a
representation of the analytic extension $(\hat \ar ,\hat \dar
)$ of $(\hat A, \hat \D)$ on the space~$H$. There is a unique
$*$-homomorphism $\hat{\pi} :\hat{A}\rightarrow B(H)$ such
that  $\hat{\pi} (a )\hat{\Lambda}(b )=\hat{\Lambda}(a b)$,
for all $a,b\in\hat{A}$. Moreover, $\hat\pi$ is faithful and
non-degenerate. Let $\hat A_{\mathrm{r}}$ be the norm closure
of $\hat{\pi} (A)$ in $B(H)$, so $\hat A_{\mathrm{r}}$ is a
non-degenerate C*-subalgebra of $B(H)$. One can show that
$W\in M(B_0(H)\otimes \hat \ar )$ and that $\hat
A_{\mathrm{r}}$ is the norm closure of the linear space ${
\{(\omega\otimes \iota)(W)\mid \omega\in B_0(H )^*\}}$. Define
a linear map $\hat \D_{\mathrm{r}} :\hat
A_{\mathrm{r}}\rightarrow M(\hat A_{\mathrm{r}}\otimes\hat
A_{\mathrm{r}} )$ by setting $\hat \D_{\mathrm{r}} (a )=W
(a\otimes 1)W^*$, for all ${a \in \hat A_{\mathrm{r}}}$. Then
$\hat \D_{\mathrm{r}}$ is  the unique $*$-homomorphism ${\hat
\D_{\mathrm{r}} : \hat A_{\mathrm{r}} \rightarrow M(\hat
A_{\mathrm{r}} \otimes\hat A_{\mathrm{r}})}$ such that, for
all $a \in \hat{A}$ and ${x \in \hat{A} \otimes \hat{A}}$,
\[\hat \D_{\mathrm{r}} (\hat\pi(a)) (\hat\pi \otimes \hat\pi) (x) =
(\hat\pi \otimes \hat\pi) (\hat{\D} (a)x)\] and
\[(\hat\pi \otimes \hat\pi) (x)\hat \D_{\mathrm{r}} (\hat\pi(a))
= (\hat\pi \otimes \hat\pi) (x \hat{\D}(a)).\] Note that one
can show that $W\in M(\ar \otimes \hat \ar )$ and ${(\dar
\otimes \iota)(W)}=W_{13}W_{23}$.

\vspace{2ex} An algebraic quantum group $(A, \Delta)$ is of
{\em compact type} if $A$ is unital, and of {\em discrete
type} if there exists a non-zero element $h \in A$ satisfying
$ah = ha = \varepsilon(a)h$, for all $a \in A$.
It is known that $(A, \Delta)$ is of compact type
(respectively, of discrete type) if, and only if, its dual
$(\hat A, \hat \Delta)$ is of discrete type (respectively, of
compact type).

\vspace{2ex}
We use the symbol $\hat{M}$ to denote the von~Neumann algebra
generated by~$\hat A_{\mathrm{r}}$, so that $\hat
A_{\mathrm{r}}$ and $\hat\pi(\hat A)$ are weakly dense in
$\hat{M}$. As with $\dar$, since $\hat \D_{\mathrm{r}}$ is
unitarily implemented, it has a unique extension to a weakly
continuous unital $*$-homomorphism ${\hat \D_{\mathrm{r}}
:\hat{M}\rightarrow \hat{M}\bar{\otimes}\hat{M}}$, given
explicitly by $\hat\D_{\mathrm{r}}(a)={W(a\otimes 1)W^*}$, for
all $a\in \hat M$.

\vspace{2ex}
It should be noted that both $M$ and $\hat{M}$ are in the
standard representation. This follows easily from \cite{KD}
and standard von~Neumann algebra theory (see~\cite{SZ}, for
example). As a consequence, all normal states on these algebras
are (restriction of) vector states.

\vspace{3ex}
We now recall the definition of co-amenability of an \aqg. Suppose
that $\AA$ is an \aqg\
and let $(H,\Lambda)$ be a GNS~pair associated to a left Haar
integral. As the representation ${\pi\colon A\to B(H)}$ is
injective, we can use it  to endow $A$ with a C*-norm by
setting $\n a=\n{\pi(a)}$, for $a\in A$. Following \cite{BMT2} (see
also \cite{BMT}), we say that $(A, \D)$ is {\em co-amenable}
if its co-unit $\e$ is norm-bounded
with respect to this norm. Several characterizations of
co-amenability are given in  \cite{BMT2}.
We just remind that the \aqg of compact type associated to the group algebra
of a
discrete group $\Gamma$ is co-amenable according to this
definition if, and only if, $\Gamma$ is amenable.
Also, co-amenability is automatic in the case of
a discrete type \aqg.

\vspace{2ex}
We also recall from \cite{BMT2} the definition of amenability for an \aqg.
Let $\AA$ be an \aqg\ with von~Neumann algebra~$M.$ \\ A {\em
right-invariant mean} for $(M, \D)$ is a state $m$ on $M$ such
that
\[m((\iota\bar{\otimes}\omega) \D(a))=\omega (1)m(a),\]
for all $a\in M$ and $\omega\in M_*$. A {\em left-invariant
mean} is defined analogously. We say that $\AA$ is {\em
amenable} if $(M, \D)$ admits a right-invariant mean. Using
the existence of the anti-unitary antipode $R$ on $(M, \D)$
(\cite{KD,KV}), this is easily seen to be equivalent to
requiring that $(M, \D)$ admits a left-invariant mean.
The  \aqg\ associated to the algebra of complex functions with finite
support on a
discrete group $\Gamma$ is amenable if, and only if, the group $\Gamma$ is
amenable, by the very definition of the amenability of a group.
Amenability is automatic for an \aqg\ $\AA$ of compact type.

\vspace{2ex}

We end this section with some technical lemmas.

\vspace{2ex}

We denote by $\delta$ the modular ``function'' of $\AA.$ Especially,
$\d$ is an invertible, self-adjoint element of $ M(A)$ satisfying

\[ \D(\delta) = \delta \otimes \delta, \  \e(\delta) =1, \  S(\delta) =
\delta^{-1} .\]

Further, there exists $\mu \in \bT $ such that

\[ \varphi S (a) = \varphi( a \delta ) = \mu \varphi( \delta a ) \]

for all $a \in A.$

\begin{lem} \label{techlem1}
    Let $a \in A.$ Then we have

    \vspace{2ex}

\hspace{5ex}    (1) $\hat{a}^{*} = (S(a)^* \delta )^{\wedge} ,$

    \vspace{2ex}

\hspace{5ex}    (2) $\hat{a} \delta^{-1} = \mu (a \delta^{-1})^{\wedge},$

    \vspace{2ex}

\hspace{5ex}    (3) $ \widehat{S(a)} \rho^{-1} S = \widehat{a \delta},$

    \vspace{2ex}

\hspace{5ex}    (4) $ \widehat{\rho(a^*)} S^{-1} = \mu^{-1} (S(a^{*})
    \delta)^{\wedge}. $

    \end{lem}

    \demo

    \vspace{2ex}
 
 (1) See \cite[Lemma 7.14]{KD}.
 
 \vspace{2ex}
 
 (2) Observe that
 \[ \varphi(\delta^{-1} b) =  \varphi(\delta^{-1} b \delta^{-1} \delta) =
 \mu \varphi(\delta  \delta^{-1} b \delta^{-1}) = \mu \varphi(b
 \delta^{-1} )\]
 
 for all $b \in A.$ Hence we get
 
 \[ (\hat{a} \delta^{-1}) (c) = \hat{a} ( \delta^{-1} c) = \varphi (
 \delta^{-1} c a ) \]
 \[ = \mu \varphi ( c a \delta^{-1} ) = \mu ( a \delta^{-1} )^{\wedge}
 (c) \]
 for all $ c \in A.$
 
 \vspace{2ex}
 
 (3) We have
 \[ (\widehat{S(a)} \rho^{-1} S ) (b) = \varphi ( \rho^{-1}( S(b)) S(a)) =
    \varphi ( S(a) S(b) ) \]
  \[= \varphi S (b a) = \varphi ( b a \delta) = \widehat{a \delta} (b) \]
  for all $b \in A.$
 
  \vspace{2ex}
 
  (4) We have
 
  \[  ( \widehat{\rho(a^*)} S^{-1}) (b) = \varphi ( S^{-1}(b) \rho(a^*) )
      = \varphi ( a^* S^{-1}(b) ) \]
   \[ = \varphi S ( a^* S^{-1}(b) \delta^{-1} ) = \varphi ( \delta b S(a^*))
      = \mu^{-1} \varphi ( b S(a^* \delta) \]
      \[ = \mu^{-1} (S(a^{*}) \delta)^{\wedge} (b)\]
      for all $b \in A.$~\qed

\vspace{2ex} 

 Let now $\hat \rho$ denote the automorphism of $\hat A$ satisfying
 \[ \hat \psi ( \hat a \hat b ) = \hat \psi ( \hat b \hat{\rho}( \hat
 a)  ) \]
 for all $a,b \in A.$ (The existence of $\hat \rho$ is proved in a
 similar way as the existence of $\rho.$)

 \begin{lem} \label{techlem2}
     Let $\a \in A.$ Then we have

     \vspace{2ex}

\hspace{5ex}     (1) $\hat \rho ( \hat a) = ( S^{2}(a) \delta^{-1}
)^{\wedge} $

\vspace{2ex}

\hspace{5ex}     (2) $ \hat \rho ( \widehat{ S(a^*)}^* ) = \hat a $

     \end{lem}

 \demo
 \vspace{1ex}

 (1) For all $b \in A$ we have

 \[ \hat \psi ( \hat{ b}^* \hat \rho (\hat a)) =
     \hat \psi ( \hat a \hat{ b}^* ) =
     \hat \psi ( (\hat{ a}^*)^* \hat{b}^* ) \]
      \[ = \hat \psi ( ( (S(a)^* \delta) ^{\wedge})^{*} (S(b)^*
     \delta)^{\wedge}) \ \ \ \ (\textnormal{ using Lemma \ref{techlem1},
     (1)}) \]
     \[ = \f ( (S(a)^* \d)^* S(b)^{*}\d ) = \f ( \d^* S(a) S(b)^*
     \d) = \f ( \d S(a) S^{-1}(b^*) \d)\]
     \[ = \f S ( \d S(a) S^{-1}(b^*)) = \f ( b^* S^2(a) \d^{-1}) =
     \hat \psi  (\hat {b}^* (S^2(a)\d^{-1})^{\wedge}), \]

     and the assertion follows from the faithfulness of $\hat \psi.$

\vspace{2ex}

(2) Observe that
\[ S^2( S^2(a^*)^*) = S( S( S^2(a^{*})^{*})) =
S(S^{-1}(S^2(a^*))^*) \]
\[ = S(S(a^{*})^*) = a .\]
Hence, using (1) from Lemma \ref{techlem1}, and (1) above, we get

\[ \hat \rho(\hat{S(a^{*})}^*)
    = \hat \rho ( (S(S(a^*))^* \d)^{\wedge} ) = \hat \rho (
    ( S^2(a^{*})^* \d)^{\wedge} ) \]
     \[ = ( S^2 (S^2(a^*)^* \d) \d^{-1} )^{\wedge} )
    = (S^2( S^2(a^*)^*) S^2(\d) \d^{-1})^{\wedge}
    = (S^2 (S^2(a^{*})^*) )^{\wedge} = \hat a\, \]

    as desired.\qed

    \vspace{2ex}
 
    \begin{lem} \label{techlem3}
 
 Define $F: \hat A \to \hat A $ by $ F(\hat a) = \widehat{S(a^*)}.$
 Then $F$ is antilinear, antimultiplicative and involutive.
 
 \end{lem}
 
 \demo We only show antimultiplicativity as the other two properties are
 easily checked.
 \vspace{1ex}
 
 Let $a,b \in A.$ Write $a \otimes b = \sum_{i=1}^{n} \D(p_{i})
 (q_{i} \otimes 1) $ for some $p_{i}, q_{i} \in A,i = 1 \ldots n.$ Then
 \[ (\hat a \hat b) =  \sum_{i=1}^{n} \f (q_{i}) \hat{p_{i}}. \]
 
 Indeed, we have
 
 \[ (\hat a \hat b) (c) = (\hat a \otimes \hat b) \D(c) = (\f \otimes
 \f ) ( \D(c) (a \otimes b) )\]
 \[ = (\f \otimes \f ) ( \D(c)  \sum_{i=1}^{n} \D(p_{i}) (q_{i} \otimes 1))
 = \sum_{i=1}^{n}(\f \otimes \f ) ( \D(c p_{i})
 (q_{i} \otimes 1)) \]
 \[ = \sum_{i=1}^{n} \f(c p_{i}) \f(q_{i})
  = \f ( c (\sum_{i=1}^{n} \f (q_{i}) p_{i}) = ( \sum_{i=1}^{n} \f
 (q_{i}) \hat{p_{i}} ) (c) \]
 
 for all $c \in A.$
 Using this expression, we get
 
 \[ F( \hat a \hat b) (c) = ( \sum_{i=1}^{n} F(\f
 (q_{i}) \hat{p_{i}}) ) (c) \]
 \[ = \sum_{i=1}^{n} \f (q_{i}^*) \widehat{S(p_{i}^*} ) (c)
    =  \sum_{i=1}^{n} \f (q_{i}^*) \f (c S(p_{i}^*) ) , \]
 
    while
 
    \[( F(\hat b ) F(\hat a)) (c) = (\widehat{S(b^*)}
    \widehat{S(a^*)}) (c) = (\widehat{S(b^*)} \otimes
    \widehat{S(a^*)}) \D(c)  \]
     \[ = (\f \otimes \f) ( \D(c) (S(b^{*}) \otimes S(a^*)) )
        = \sum_{i=1}^{n}  (\f \otimes \f)
        ( \D(c) (S \otimes S) ( (1 \otimes q_{i}^*) \chi \D
        (p_{i}^{*}))) \]
    \[ = \sum_{i=1}^{n}  (\f \otimes \f)
        ( \D(c) ((S \otimes S) \chi \D
        (p_{i}^{*})) ( 1 \otimes S(q_{i}^*)) ) \]
    \[  =   \sum_{i=1}^{n}  (\f \otimes \f)
        ( \D(c) (\chi (S \otimes S)  \D
        (p_{i}^{*})) ( 1 \otimes S(q_{i}^*)) )
        \]
    \[  = \sum_{i=1}^{n}  (\f \otimes \f)
        ( \D(c S(p_{i}^*))  ( 1 \otimes S(q_{i}^*)) )
        = \sum_{i=1}^{n}  (\f S \otimes \f)
        ( \D(c S(p_{i}^*))  ( \d^{-1} \otimes S(q_{i}^*)) ) \]
    \[ = \sum_{i=1}^{n}  (\f S \otimes \f)
        ( \D(c S(p_{i}^*) \d^{-1})  ( 1 \otimes  \d S(q_{i}^*)) )
      = \sum_{i=1}^{n}  \f S (c S(p_{i}^*) \d^{-1}) \f (\d
      S(q_{i}^*)) \]
    \[ =\sum_{i=1}^{n}  \f (c S(p_{i}^*)) \f S (q_{i}^* \d^{-1})
    = \sum_{i=1}^{n}   \f  (q_{i}^* ) \f (c S(p_{i}^*)) \]
 
    for all $c \in A,$ and the antimultiplicativity of $F$ follows.~\qed
 
\vspace{2ex}

\begin{lem} \label{techlem4}
    Let $a \in A.$ Pick $c \in A$ such that $\widehat{S(a^{*})} =
   \hat c \widehat{S(a^{*})} .$ Then we have
    \[ (\iota \otimes \f) ( (1 \otimes c^* ) \D(a) ) = a \]
 \end{lem}

 \demo Using Lemma \ref{techlem3}, we get $\hat a \widehat{ S(c^*)} = \hat
a.$
 Now, to prove the assertion, it is clearly enough  to show that

 \[ (\f \otimes \f)( ( b^* \otimes c^*) \D(a) ) = \f (b^* a) \]
holds for all $b \in A.$
 This may be established as follows :
 \vspace{2ex}

 $ (\f \otimes \f)( ( b^* \otimes c^*) \D(a) ) = \f ( b^* (\iota
 \otimes \f c^*) \D(a) ) $

 \[ = \f S (b^* ( ( \iota \otimes \f c^*) \D(a) ) \delta^{-1}) = \f (
 \delta (( S \otimes \f c^*) \D(a)) S(b^*)) \]
 \[ = \f ( \delta (( S \otimes \f c^* S^{-1}S) \D(a)) S(b^*)) =
 \f ( \delta ((   \f c^* S^{-1} \otimes \iota) \D(S(a))) S(b^*)) \]
 \[ = \f ( ((   \f c^* S^{-1} \otimes \iota) \D(S(a))) S(b^*) \rho(\delta))
 = \mu^{-1} \f ( ((   \f c^* S^{-1} \otimes \iota) \D(S(a)))
 S(\delta^{-1} b^*) ) \]
 \[ = \mu^{-1} \f ( \rho^{-1} S (\delta^{-1} b^*) (( \f c^* S^{-1}
 \otimes \iota) \D(S(a)) ) ) \]
\[ = \mu^{-1} \f ( ( ( \f c^* S^{-1} S \otimes \iota )
 \D(\rho^{-1}S(\delta^{-1}b^*))) S(a) ) \]
 ( where we have used strong left invariance of $\f$)
 \[ = \mu^{-1}  ( \f c^* \otimes S(a) \f) \D \rho^{-1} S ( \delta^{-1}
b^*)  = \mu^{-1}  ( \widehat{\rho(c^*)} \otimes \widehat{S(a)})  \D
\rho^{-1} S ( \delta^{-1}
b^*)  \]
\[ = \mu^{-1}  ( \widehat{\rho(c^*)} \otimes \widehat{S(a)}) (S^{-2}
\otimes \rho^{-1}) \D  S ( \delta^{-1} b^*)  \]
\[ =  \mu^{-1}  (\widehat{S(a)}  \otimes \widehat{\rho(c^*)}) (S^{-2}
\otimes \rho^{-1})  (S \otimes S) \D ( \delta^{-1} b^*)  \]
\[ =  \mu^{-1}  ( \widehat{S(a)} \rho^{-1} S \otimes \widehat{\rho(c^*)}
S^{-1})  \D ( \delta^{-1} b^*)  =   \mu^{-2} ( \widehat{a \delta}
\otimes \widehat{S(c^*)\delta})  \D ( \delta^{-1} b^*) \]
(where we have used Lemma \ref{techlem1}, (3) and (4))
\[ =   \mu^{-2} ( \widehat{a \delta}
\otimes \widehat{S(c^*)\delta}) (\delta^{-1} \otimes \delta^{-1}) \D (b^*)
=   \mu^{-2} \mu^{2} ( \hat{a} \otimes \widehat{S(c^*)})  \D (b^*) \]
(where we have used Lemma \ref{techlem1}, (2))

\vspace{2ex}

 $ =(\hat a \, \widehat{S(c^*)}) (b^*) = \hat a (b^*) = \f( b^* a), $
 where we have used that $(\hat a  \, \widehat{S(c^*)}) = \hat a .$

 \vspace{2ex}

 This finishes the proof.~\qed

  \begin{lem} \label{techlem5} Let $a, b \in A.$  Pick $\ a_{1},b_{1},
  \ldots, a_{n}, b_{n} \in A$ such that
     $$ \D(b) (a \otimes 1) = \sum_{k=1}^{n} a_{k} \otimes b_{k}.$$
     Then, for all $\om \in \hat A,$ we have

      $$( (a \om) S) \hat b = \sum_{k=1}^{n} \om (a_{k}) \hat{b}_{k}.$$

\end{lem}
 
\demo Let $c \in A, \om \in \hat A.$ Then
$$ (((a \om)S) \hat b) (c) = (((a \om) S) \odot (b \f)) \D(c)= (
a\om)S  ( \iota \odot (b \f)) \D(c)$$
$$=(a \om )S(\iota \odot \f) ( \D(c) (1 \otimes b))=
(a \om ) ( \iota \odot \f)((1 \otimes c) \D(b)),$$
using  strong left invariance of $\f$.
Hence,
$$ (((a \om)S) \hat b) (c) = (\om \odot \f)((1 \otimes c) \D(b) (a
\otimes 1)) = (\om \odot \f)((1 \otimes c)(\sum_{k=1}^{n}
a_{k} \otimes b_{k}))$$
$$ = \sum_{k=1}^{n} \om (a_{k}) \f (c b_{k}) = (\sum_{k=1}^{n} \om (a_{k})
  \hat{b}_{k}) (c).$$
  This shows the assertion.~\qed

 \begin{lem} \label{techlem6}  Let $a, b, c \in A.$ Write
     $$ \D(b) (a \otimes 1) = \sum_{k=1}^{n}
a_{k} \otimes b_{k} \ , \ b \otimes c = \sum_{i=1}^{r}
 \D(p_{i}) (q_{i}\otimes 1) $$ for some $a_{1},b_{1}, \ldots, a_{n},
 b_{n}, p_{1},q_{1}, \ldots, p_{r}, q_{r} \in A.$
 Further, for
 each $i = 1, \ldots , r, $ write
 $$\D (p_{i}) (a \otimes 1) =  \sum_{j=1}^{s(i)}
x_{ij} \otimes y_{ij}$$ for some $x_{i1},y_{i1}, \ldots, x_{is(i)},
y_{is(i)} \in A.$
Then

$$ \sum_{i=1}^{m} \sum_{j=1}^{s(i)}  \varphi (q_{i})   \hat{y}_{ij}
\otimes x_{ij} =\sum_{k=1}^{n}
  \hat{b}_{k} \hat c \otimes a_{k} \, .$$

\end{lem}
 
 \demo   Set $ X = \sum_{i=1}^{m} \sum_{j=1}^{s(i)}  \varphi (q_{i})
\hat{y}_{ij}
\otimes x_{ij}$ and $Y=\sum_{k=1}^{n}
  \hat{b}_{k} \hat c \otimes a_{k},$  which are both elements in $\hat
  A \odot A.$ To show that $X=Y$, it suffices, using separation, to
  prove that
  $(e\otimes (\f f)) X =(e\otimes (\f f)) Y$ for all $e, f \in A.$
  Note that we regard here $e$ as an element of the double dual of $A.$

  \vspace{1ex}
Let $e, f \in A.$  Then

$$(e\otimes (\f f)) X = \sum_{i,j} \f(q_{i}) \f(e y_{ij}) \f(f
x_{ij})$$
$$= \sum_{i,j} (\f \otimes \f \otimes \f) ((1 \otimes e \otimes f)
(q_{i} \otimes y_{ij} \otimes x_{ij}))$$
$$ = \sum_{i}( \f \otimes \f \otimes \f) ((1 \otimes e \otimes f)
(q_{i} \otimes (\D_{\mathrm{op}}(p_{i}) (1\otimes a))))$$
$$= (\f \odot (\f e) \odot (a \f f)) ( \sum_{i} q_{i} \otimes
\D_{\mathrm{op}}(p_{i}))$$
$$= (\f \odot (\f e) \odot (a \f f)) ((\iota \odot
\D_{\mathrm{op}}) ( \sum_{i} q_{i} \otimes p_{i}) )$$
$$= (\f \odot (\f e) \odot (a \f f)) ((\iota \odot
\D_{\mathrm{op}}) (S^{-1} \odot \iota) ( \D (c)) ( b \otimes 1)   ), $$
using here the formula established in \cite[Proposition 2.2]{KD} at the last
step.

\vspace{1ex}
Continuing this computation, we get
$$(e\otimes (\f f)) X = ( ( (b \f ) S^{-1} \odot (\f e) \odot (a \f
f)) ((\iota \odot \D_{\mathrm{op}}) \D (c))$$
$$ = ( ( (b \f ) S^{-1} \odot (a \f f)\odot (\f e)  ) ((\iota \odot \D) \D
(c))$$
$$ = ( ( (b \f ) S^{-1} \odot (a \f f) \odot (\f e)  ) (( \D \odot \iota )
\D (c))$$
$$= (( (b \f ) S^{-1} \odot (a \f f)) \D ( (\iota \odot (\f e)) \D
(c))$$
$$= (( (b \f ) S^{-1} \odot (a \f f)) \D ( (\iota \odot \f)
((1\otimes e) \D (c)))$$
$$= (( (b \f ) S^{-1} \odot (a \f f)) \D ( S( (\iota \odot \f)
(\D( e) (1 \otimes c))  )),$$
using  strong left invariance of $\f.$

\vspace{1ex}
This gives
$$(e\otimes (\f f)) X = (( (a \f f) S)  \odot (b \f )) \D ((\iota \odot \f)
(\D( e) (1 \otimes c))  ))$$
$$=(( (a \f f) S ) \odot (b \f )) \D ((\iota \odot \hat c) \D(e))$$
$$= (( (a \f f) S)  \odot \hat b \odot \hat c) ) (\D \odot \iota)
\D(e)$$
$$= (( (a \f f) S)  \odot \hat b \odot \hat c) ) (\iota \odot \D)
\D(e)$$
$$= (( (a \f f) S)  \odot (\hat b \hat c)) \D(e) =  (( (a \f f) S)
\hat b \hat c) (e).$$
On the other hand, we have
$$ (e \odot (\f f)) Y = \sum_{k=1}^{n}
  ((\hat{b}_{k} \odot \hat c) \D(e)) (\f f) ( a_{k})$$
  $$ = \sum_{k=1}^{n}
  (\hat{b}_{k}  \hat c) (e) (\f f) ( a_{k}) = ((\sum_{k=1}^{n}(\f f) (
  a_{k}))
  \hat{b}_{k})  \hat c) (e).$$
  Hence, the assertion will be proved if we can show that
  $$( (a \f f) S) \hat b = \sum_{k=1}^{n}(\f f) (a_{k}) \hat{b}_{k}$$
  holds for all $f \in A.$ But this clearly follows from Lemma
  \ref{techlem5}, and the proof is finished.~\qed
 
\section{Non-degenerate representations and unitary generators}

 We let $\AA$ be an \aqg\ throughout this section and use notation and
 terminology introduced in the previous section.

 \vspace{2ex}
 Following \cite{Ku2}, we first introduce the universal C$^{*}$-algebraic
 quantum group $(\au,\dau)$ associated to $\AA$. The
 C$^{*}$-algebra $\au$  is the completion of $A$ with respect to the
 C$^{*}$-norm $\| \cdot \|_{\mathrm{u}} $ on $A$ defined by
 \[ \| a \|_{\mathrm{u}} = \sup \{ \| \varphi(a) \| \ | \  \varphi \
 \textnormal{is a
 $*$-homomorphism from $A$ into some C$^{*}$-algebra}  \}. \]
 (The non-trivial fact that this gives a well-defined norm on $\au$ is
 shown in \cite{Ku2}.)
  Let $\pu$ denote the identity
 mapping from $A$ into $\au$. The co-product map  $\dau $
 is defined in such way that it is the unique non-degenerate
 $*$-homomorphism $\dau : \au \to M(\au \otimes \au)$ satisfying

  \[ (\pu \odot \pu)(x) \dau (\pu(a)) =
 (\pu \odot \pu)(x \Delta(a))   \]
 and
 \[ \dau  (\pu(a)) (\pu \odot \pu)(x) =
 (\pu \odot \pu)(\Delta(a) x)   \]
 for all $a \in A$ and $x \in A \otimes A.$

 \vspace{1ex}
 The universality
 of $\au$ makes it possible to extend uniquely from $A$ to $\au$ any
$*$-homomorphism  from $A$ into
 some C$^{*}$-algebra. Especially, the co-unit $\e$ of $\AA$ extends
 to a $*$-homomorphism $\e_{\mathrm{u}} : \au \to \bC $  such that
 $\e_{\mathrm{u}} \circ \pi_{\mathrm{u}} =  \e.$ One easily checks
 that $\eu$  satisfies the co-unit property, that is,
 $$(\eu \otimes \iota)
 \dau(a) = (\iota \otimes \eu) \dau(a) = a , \ \ \ a \in \au.$$
 It follows immediately from this that $\dau$ is injective.
 Also, there exists a unique
 $*$-homomorphism $\pr$ from $\au$ onto $\ar$ satisfying $\pr
 \circ \pi_{\mathrm{u}} = \pi$. By construction, we have $(\pr \otimes
 \pr) \circ \dau = \dar \circ \pr$.

 \vspace{2ex}
 We remark that one easily verifies that our definition of
 co-amenability of $\AA$ may now be rephrased as saying that
 $\e_{\mathrm{u}}$ is weakly contained in $\pr$ , that is, ker($\pr$)
$\subset$ ker($\e_{\mathrm{u}}$).

 \vspace{2ex}
 One of the off-springs of \cite{Ku2} is that there is a bijective
 correspondence between non-degenerate $*$-homomorphisms of $\au$ and
 unitary corepresentations of $(\hat \ar, \hat \dar)$ (this may be
 seen by
 combining results from section 7 and section 13 in \cite{Ku2}). Kustermans
has also
 established a similar result for more general locally compact quantum
 groups in \cite{Ku3}. For completeness, we recall the definition of a
unitary
 corepresentation. Consider a
 C$^{*}$-algebraic quantum semigroup $(B,\Gamma),$ that is, a
 C$^{*}$-algebra $B$ equiped with a co-product map $\Gamma.$ Then a
 \emph{unitary corepresentation} of $(B,\Gamma)$ on a (non-zero) Hilbert
space $K$ is
 a unitary element $U \in M(B \otimes B_{0}(K))$ such that $(\Gamma  \otimes
 \iota) \, U = U_{13} U_{23}.$

 \vspace{1ex}
 In this section we show that any non-degenerate $*$-representation
of $\au$ on a Hilbert space has a unitary ``generator'', that is, it
arises from some unitary corepresentation of $(\hat \ar, \hat
\Delta_{\mathrm{r, op}})$. By $\hdarop$ we mean throughout this paper the
co-product
on $\hat \ar$ \emph{opposite} to $\hdar.$ It may be seen as  a matter of
taste choosing
to work with the opposite co-product on the dual side.
 However, one reason for this choice is that this is the one tacitely
adopted by Kustermans and Vaes
\cite{KV} in the setting of locally compact quantum groups:
the ``dual'' of $(\ar,\dar)$  in their sense is
in fact precisely $(\hat \ar,\hdarop)$.

 \vspace{1ex}
For completeness, we recall  how $\hdarop$ is defined.
Let $ \hat \chi \colon \hat \ar \otimes \hat \ar  \to \hat \ar \otimes
\hat \ar $ denote the flip map. We also denote by the same symbol its
extension to a $*$-automorphism of $M (\hat \ar \otimes
\hat \ar).$ Then $\hdarop =  \hat \chi \circ \hdar$.

One can also describe $\hdarop$ with the help of a multiplicative unitary
$\hat{ W}$ related to $W,$ in the spirit of \cite{BS}. Indeed, let $ \Sigma$
denote the flip map
on $H \otimes H$ and set $\hat W = \Sigma W^{*} \Sigma$. Then one checks
readily that
\[ \hdarop(y) = \hat W^{*}  (1 \otimes y) \hat W \]
for all $y \in \har$. We may use this formula to extend
$\hdarop$ to a  map from  $\hat M$ into $\hat{M} \bar{\otimes} \hat M$,
which we also  denote by
the same symbol. The pair $(\hat M, \hdarop)$ is then a Hopf von
Neumann algebra.

\vspace{1ex}
We  equip $\hat{\ar}^*$ (resp. $\hat M_{*}$) with the
product induced by  $\hdarop$, that is $\om \eta = (\om \otimes \eta )
\hdarop$ (resp. $\om \eta = (\om \bar{\otimes} \eta )
\hdarop$). It is then straightforward to check that the Banach space
$\hat{\ar}^*$ (resp. $\hat M_{*}$) is a Banach algebra under this
product.

\vspace{1ex}
Our approach relies on the following fundamental result, which
takes advantage of the fact that we are dealing with algebraic quantum
groups.

\begin{pro} \label{Qmap}
    Let notation be as above. Then :
    \vspace{1ex}

    (a) There exists a (unique) injective homomorphism $Q: A \to \hat
    M_{*}$ satisfying
    \[ Q(a) [\hat \pi( \hat b)] = \varphi ( S^{-1}(a) b)  \]
    \hspace{5ex} for all $a,b\in A$.
    If $a \in A$, and $c \in A$ is chosen  so that
    \[ \hat c \ \widehat{S(a^{*})} = \widehat{S(a^{*})}, \]
      \hspace{5ex} then $Q(a) = \om_{\Lambda(a), \Lambda(c)}$ (restricted to
$\hat
    M$).

    \vspace{2ex}
    (b) The algebra $Q(A)$ is norm-dense in $\hat M_{*}$.

    \vspace{2ex}
    (c) If $\Theta$ denotes the canonical $*$-homomorphism from $\AA$
    onto its double
   \par{} \hspace{3ex} dual, then $Q(a) \circ  \hat \pi =
    \Theta(S^{-1}(a))$ for all $a \in A$.

    \vspace{2ex}
    (d) Let $Q_{\mathrm{r}}$ denote the map from $A$ into $\hat
    \ar^{*}$ obtained by restricting each $Q(a)$ \par{} \hspace{3ex}
    to $ \hat \ar$. Then
    $Q_{\mathrm{r}}$  is also an injective homomorphism.

    \end{pro}

    \vspace{1ex}
    \demo Let $a \in A$. We define a linear functional $P(a)$
    on $\hat \pi (\hat A)$ by
    \[ P(a) [\hat \pi( \hat b)] = \varphi ( S^{-1}(a) b)  \]
    for all $b \in A$.
    Now choose $c \in A$ such that $\hat c \ \widehat{S(a^{*})} =
\widehat{S(a^{*})}$.
    Then we have
    \[ P(a) [\hat \pi( \hat b)] = \om_{\Lambda(a), \Lambda(c)}(\hat \pi(
    \hat b)) \]
    for all $b \in A$. Indeed, let $b \in A$. Then
    \[ P(a) [\hat \pi( \hat b)] = \varphi ( S^{-1}(a) b) = \varphi (
    (S^{-1}(a)^{*})^{*} \  b) = \varphi( S(a^{*})^{*} \ b) \]
 \[= \hat \psi ( (\widehat{S(a^{*})})^{*} \ \hat{b}) =
        \hat \psi ( (\widehat{S(a^{*})})^{*} \hat{c}^{*} \ \hat{b})
 = \hat \psi (  \hat{c}^{*} \ \hat{b} \ \hat \rho (
 (\widehat{S(a^{*})})^{*}) ) \]
\[ = \hat \psi (  \hat{c}^{*} \ \hat{b} \ \hat a ) \ \ \ \ \ \
 \textnormal{( using Lemma \ref{techlem2}) } \]
 \[ = ( \ \hat \Lambda( \hat b \hat a) , \hat \Lambda( \hat c) \ )
    = ( \  \hat \pi (\hat b) \hat \Lambda ( \hat a) , \hat \Lambda ( \hat 
    c) \ )
    = ( \  \hat \pi (\hat b)  \Lambda( a ) , \Lambda( c ) \ ) \]
  \[= \om_{\Lambda(a), \Lambda(c)}(\hat \pi(\hat b)), \]
   as asserted.
 
   It follows clearly from the formula just established that $P(a)$
   has a unique extension $Q_{\mathrm{r}}(a) \in \hat \ar$ and also a
   unique extension to $Q(a) \in \hat M_{*}$, both determined by
   restricting suitably $\om_{\Lambda(a), \Lambda(c)}$.
 
\vspace{2ex}
 We show now that assertion (c) holds. Let $a \in A$. Then we have

 \[ (Q(a) \circ \hat \pi) ( \hat b) = Q(a) ( \hat \pi ( \hat b)) =
 \varphi (S^{-1}(a) \ b) \]
 \[ = \hat b ( S^{-1}(a) ) = (\Theta (S^{-1}(a))) (\hat b) \]
 for all $b \in A$. Hence, $ Q(a) \circ \hat \pi = \Theta(
 S^{-1}(a))$, as desired.

\vspace{2ex}
The map $Q : a \to Q(a) $ from $A$ into $\hat M_{*}$ is clearly
linear. We show  that $Q$ is multiplicative.

Let $a_{1}, a_{2} \in A$. For all $b \in A$ we have

\[ ( Q(a_{1})  \ Q(a_{2}) ) ( \hat \pi ( \hat b) ) =
   ( Q(a_{1})  \bar{\otimes} Q(a_{2}) ) ( \hdarop(\hat \pi ( \hat b)) )
   \]
 \[ = ( Q(a_{2})   \bar{\otimes} Q(a_{1}) ) ( \hdar (\hat \pi ( \hat
   b) )) = ( Q(a_{2})   \bar{\otimes} Q(a_{1}) ) ((\hat \pi \odot \hat
   \pi) \hat \D (\hat b) ) \]
 \[ = ( \Theta ( S^{-1}(a_{2})) \odot \Theta ( S^{-1}(a_{1})) ) \hat \D
(\hat b)
      = \hat \D (\hat b) (S^{-1}(a_{2}) \otimes S^{-1}(a_{1})) \]
 \[ = \hat b ( S^{-1}(a_{2}) \ S^{-1}(a_{1}) ) = \hat b (
      S^{-1}(a_{1}a_{2}) ) \]
 \[ = \varphi( S^{-1}(a_{1}a_{2}) \ b) = ( Q(a_{1}a_{2})) (\hat \pi (
\hat b) ).\]

      As $\hat \pi (\hat A)$ is weakly dense in $\hat M,$ the
      multiplicativity of $Q$ follows.

      \vspace{2ex}
      To finish the proof of (a), it remains only to show that $Q$ is
      injective. Let $a \in A$ and assume that $Q(a) =0$.

      Then, for all $b \in A$, we have  $0 = Q(a) (\hat \pi (\hat
      b)) = \varphi ( S^{-1}(a) b) $. Hence, inserting $b =
      (S^{-1}(a))^{*}$, we get $\varphi(S^{-1}(a) \ (S^{-1}(a))^{*})= 0$,
      that is,  $S^{-1}(a)=0$ since $\varphi$ is faithful on $A$. Thus,
      $a=0$, as desired.

      \vspace{2ex}
      A little thought shows that assertion (d) is also established by
      the arguments given so far. We finally prove assertion (b).

       \vspace{2ex}
      We first show that $Q(A)= \{ {\om_{\Lambda(e), \Lambda(f)}}_{|
      \hat M} \  |  \ e, f \in A \}.$

      The inclusion $\subset$ is obvious from what we already have
      seen. To prove the reverse inclusion, let $e, f \in A$.

      Then, choosing $d \in A$ such that $ \widehat{ (S^{-1}(d))^{*}} =
      \hat f \hat \rho^{-1}(\hat e)^{*},$ we have
     \[ {\om_{\Lambda(e), \Lambda(f)}}_{| \hat M} = Q(d).\] Indeed, for
     all $b \in A$, we have
     \[ \om_{\Lambda(e), \Lambda(f)}( \hat \pi ( \hat b)) = ( \hat \pi (
     \hat b) \ \Lambda(e), \Lambda(f) ) = ( \hat \pi (
     \hat b) \ \hat \Lambda(\hat e), \hat \Lambda( \hat f) ) \]

     \[ = ( \hat \Lambda( \hat b \ \hat e), \hat \Lambda( \hat f) ) =
     \hat \psi ( \hat f ^{*} \ \hat b \ \hat e) = \hat \psi ( \hat
     \rho^{-1}(\hat e) \ \hat f ^{*} \ \hat b)\]

     \[= \hat \psi ( (\widehat{ S^{-1}(d)^{*}})^{*} \ \hat b) = \varphi (
     S^{-1}(d) b) = Q(d) (\hat \pi (\hat b) ). \]
     This proves the reverse inclusion. Now,
     since $\hat M$ acts standardly on $H$, we have $\hat M_{*} = \{
     {\om_{u,v}}_{| \hat M} \ | \ u,v \in H\} .$ Further, the
     following inequality is probably well-known (and easy to prove):

     Let $u,v \in H, a,b \in A$. Then

     \[ \|{\om_{u,v}}_{| \hat M} - {\om_{\Lambda(a),\Lambda(b)}}_{| \hat
     M} \| \leq \| u - \Lambda(a) \| \  \| v \| + \| \Lambda(a) \| \
     \| v - \Lambda(b) \| \]
      As $\Lambda(A)$ is dense
     in $H$, the norm-density of $Q(A)$ in $\hat M_{*}$  clearly
     follows.
     This finishes the proof of (b), and thereby of the
     proposition.~\qed

  \begin{thm}  \label{generator}

      Let $U$ be a unitary corepresentation of $(\hat \ar,\hdarop)$
      on a Hilbert space $K.$ Define $\pi_{U} : A \to B(K) $ by
\[\pi_U(a) = (\qr(a) \otimes \iota) \, U, \  a \in A. \]

Then $\pi_{U}$ is a non-degenerate  $*$-representation of $A$ on $K,$
that is, it is a $*$-homomorphism of $A$ into $B(K)$ which is non-degenerate
in
the sense that $[\pi_{U}(A) K ] = K.$

We shall also denote by $\pi_{U}$ the
associated non-degenerate $*$-representation of $\au$ on $K,$ and
call $U$  the generator of $\pi_{U}.$

\end{thm}

\demo We write $\pi$ instead of $\pi_{U}$ in this proof.

 \vspace{1ex}

Let $a \in A.$ Since $\qr(a) \in \har^{*}$, it is clear that $\pi(a)
\in B(K).$ The linearity of $\pi$ is evident. The multiplicativity
of $\pi$ follows from the corepresentation property of $U$ and the
multiplicativity of $\qr.$ Indeed, we have

\[ \pi(a b) = (\qr(a b) \otimes \iota) \, U
    = ((\qr(a) \qr(b)) \otimes \iota) \, U \]
\[ = ( (\qr(a) \otimes \qr(b)) \hdarop  \otimes \iota) \, U  =
    (\qr(a) \otimes \qr(b) \otimes \iota) (\hdarop \otimes \iota ) \
    U \]
\[    = (\qr(a) \otimes \qr(b) \otimes \iota) \, U_{13} U_{23}
 =  ((\qr(a) \otimes \iota) \, U) ((\qr(b) \otimes \iota) \, U) =
 \pi(a) \pi(b) \]
 for all $a,b \in A.$

 \vspace{2ex}

 To prove that $\pi$ is $*$-preserving, we have to adapt some
 arguments from \cite{KV}.

\noindent  We set $\hat W = \Sigma W^* \Sigma.$  As pointed out before, we have
$$ \drop (y) = \hat{W}^* (1 \otimes y) \hat W, \ \ y \in \har. $$
It follows that
$$ U_{13} \, U_{23} = (\drop \otimes \iota) U = \hat{W}_{12}^* \,
 U_{23} \, \hat{W}_{12}, $$
 hence that
$$(*) \ \ \ \ \ \hat{W}_{12} \, U_{13} = U_{23} \, \hat{W}_{12} \,
U_{23}^* . $$

\vspace{1ex} \noindent 
Let $ \rho \in  B_{0}(H)^* ,\eta \in B_{0}(K)^*.$ Then define $ \gamma \in
B_{0}(H)^*$
by
$$ \gamma(x) = ( \rho \otimes \eta ) ( U \, (x \otimes 1) \,
U^* ), \ \ \ x \in B_{0}(H).$$

Applying $\iota \otimes \rho \otimes \eta$ to $(*)$ above, we get
$$ ( ( \iota \otimes \rho) \,  \hat{W}) ( (\iota \otimes  \eta ) \,
U ) = (\iota \otimes \gamma) \, \hat{W} .$$
As in \cite{KV}, we can then conclude that
$ ( ( \iota \otimes \rho) \,  \hat{W}) ( (\iota \otimes  \eta ) \,
U ) \, \in D(\hs),$ and
$$\hs ( ( ( \iota \otimes \rho) \,  \hat{W}) ( (\iota \otimes  \eta ) \,
U )) = (\iota \otimes \gamma) (\hat{W}^*)$$
$$ = (\iota \otimes \rho \otimes \eta)(  U_{23} \, \hat{W}_{12}^* \,
U_{23}^*)$$
$$ = (\iota \otimes \rho \otimes \eta) (  U_{13}^* \,
\hat{W}_{12}^*) \ \ \ \ \ \ \  \textnormal{(using $(*)$)} $$
$$ = (\iota \otimes \eta) ( U^*) (\iota \otimes \rho) ( \hat{W}^* )$$
$$ = (\iota \otimes \eta) ( U^*) \hs ( (\iota \otimes \rho) \,
\hat{W} ).$$
Now, as the set $\{ (\iota \otimes \rho) \, \hat{W} \, | \, \rho \in
B_{0}(K)^* \, \}$ is a core for $\hs$, and $\hs$ is closed, this
implies that $ x( ( \iota \otimes \eta) \, U) \in D(\hs)$ and
$$ \hs( x( ( \iota \otimes \eta) \, U)) = (\iota \otimes \eta) ( U^*) \hs
(x), \ \ \ \forall x \in D(\hs).$$

>From this, we can conclude that $( \iota \otimes \eta) \, U \in
D(\overline{\hs})$ and
$$\overline{\hs}( ( \iota \otimes \eta) \, U) = (\iota \otimes \eta) (
U^*),$$
(see \cite[Remark 5.44]{KV}).

\vspace{2ex}
Let $a \in A.$ We define $Q_{r}(a)^*$ to be the linear functional on
$\hat{\pi}_{r}(\hat A)$ given by
$$Q_{r}(a)^* (\hat{\pi}_{r}(\hat b)) = \overline{Q_{r}(a) (
\hat{\pi}_{r}(\hS(\hat b)^*))}, \ \ \ b \in A.$$
Then we have

$$Q_{r}(a)^* (\hat{\pi}_{r}(\hat b)) = \overline{ (\hS(\hat b)^*)
(S^{-1}(a))}$$
$$ = (\hS( \hat b)) ( S(S^{-1}(a))^*) = (\hS( \hat b)) (a^*) =
(\hat{S}^{-1}(\hat b))(a^*)$$
$$ = \hat{b} ( S^{-1}(a^*)) = \f(S^{-1}(a^*)b)$$
$$= Q_{r}(a^*)(\hat{\pi}_{r}(\hat b)) $$
for all $b \in A.$
This shows that $Q_{r}(a)^*$ extends to an element of $\har^*$,
which is in fact equal to $Q_{r}(a^*)$.

\vspace{1ex}
Now, let  $\overline{Q_{r}(a)}$ have its usual meaning, that is
$\overline{Q_{r}(a)} \in \har^*$ is defined by
$$ \overline{Q_{r}(a) }(y) = \overline{Q_{r}(a) (y^*)}, \ \ \ y
\in \har.$$
 We have then
 $$Q_{r}(a)^*(x) = \overline{Q_{r}(a) }( \hs(x)), \ \ \ x \in D(\hs).$$
 Since $D(\hs)$ is a strict bounded core for $\hs$ (see \cite[Remark
 5.44]{KV} again), we get
 $$ Q_{r}(a)^*(x) = \overline{Q_{r}(a) }( \overline{\hs}(x)), \ \ \ x \in
 D(\overline{\hs}).$$

 Combining this with what we have seen previously, we obtain

 $$\eta( (Q_{r}(a)^* \otimes \iota) \, U) = Q_{r}(a)^*( (\iota
 \otimes \eta) \, U) $$
 $$= \overline{Q_{r}(a)}( \overline{\hs} ( ( \iota
 \otimes \eta) \, U)) = \overline{Q_{r}(a)}( ( \iota \otimes \eta) (
 U^*))$$
 $$= \eta( (\overline{Q_{r}(a)} \otimes \iota) (U^*)) = \eta (
 ((Q_{r}(a) \otimes \iota) \, U)^*),$$
 hence
 $$ \eta(\pi(a^*)) = \eta ( ( Q_{r}(a^*) \otimes \iota) \, U) $$
 $$= \eta ( (Q_{r}(a)^* \otimes \iota) \, U) = \eta (
 ((Q_{r}(a) \otimes \iota) \, U)^*)$$
 $$= \eta( \pi(a)^*).$$
 As this holds for all $ \eta \in B_{0}(K)^*,$ we can conclude that
 $\pi(a^{*}) = \pi(a)^*,$ that is,  $\pi$ is $*$-preserving.

 \vspace{2ex} Finally, we prove that $\pi $ is non-degenerate.
  Let $v\in K$ be such that $\pi(a) v = 0 $ for all $a \in A.$
  Using that $\pi$ is $*$-preserving, it is then enough to prove that $v=0.$

\noindent  Let $a,b \in A, w \in K.$ Set
 $$L = \ip{U (\LL(a) \otimes v)}{\LL(b)
 \otimes w} = \ip{ ( (\om_{\LL(a), \LL(b)} \otimes \iota) \,U)
 v} {w}.$$
 Now, it follows from the proof of Proposition \ref{Qmap} that we may
 pick $d \in A$ such that $\qr(d) $ is equal to the restriction of
 $\om_{\LL(a) , \LL(b)}$ to $\har.$  Hence, we get
 \[ L = \ip{( (\qr(d) \otimes \iota) \, U) \, v} {w} = \ip{\pi(d) \, v}{w} =
0 .\]
As this holds for all $ b \in A, w \in K,$ this implies that $U (\LL(a)
\otimes
v) = 0$ for all $a \in A.$
Thus, $\LL(a) \otimes v = U^* U (\LL(a) \otimes v) = 0$ for all $a
\in A.$ Since $\LL(A)$ is dense in  $H$, this implies $v=0,$ as
desired.~\qed

\vspace{2ex}
\emph{Remark} Let $U \in \hat M \bar{\otimes} B(K)$ be a unitary such
that $(\hdarop \bar{\otimes} \iota) \, U = U_{13} U_{23}.$ Then the
map
$\tilde{\pi}_{U}: A \to B(K)$ defined by $\tilde{\pi}_{U}(a) =
(Q(a) \bar{\otimes} \iota) \, U , a \in A,$ may also be seen to be  a
non-degenerate
$*$-homomorphism, by a similar proof. This $*$-homomorphism  extends
by universality to a
$*$-representation of $\au$ on $K.$ It will follow from our next result
that we in fact have $U \in M(\har \otimes B_{0}(K)).$ This means that $U$
is indeed a
unitary corepresentation of $(\har,\hdarop)$ and
$\tilde{\pi}_{U}= \pi_{U}.$

\vspace{2ex}
We now show that every  non-degenerate $*$-representation of
$\au$ has a unitary generator. Alternatively, one may
formulate this result for non-degenerate $*$-representations of $A.$

\begin{thm} \label{generator2}
    Let $\phi$ be a non-degenerate $*$-representation of $\au$ on some
    Hilbert space $K.$ Set $A_{\phi} = \phi (\au) =\overline{\phi(A)}
\subset B(K).$
    Then there exists a unique unitary corepresentation $U = U(\phi)$ of
    $(\har,\hdarop)$ on $K$ such that
    \[ \phi(a) = (\qr(a) \otimes \iota) \, U , \ \ \ a \in A. \]
    Hence, in the notation of Theorem \ref{generator}, we have
$\pi_{U(\phi)} = \phi .$

    \hspace{-3ex}The norm-closure of
$ \{ (\om \otimes \iota )  \, U \, | \, \om \in B_{0}(K)^* \} $ is equal to
$A_{\phi}$
and $U$ belongs to $ M(\har \otimes A_{\phi}).$

    \hspace{-3ex}Finally, we have $U(\pr)= \hat W$ and $U(\eu) = 1_{H}
\otimes 1,$
    which may be equvalently written as
    $ \pi_{\hat W} = \pr $ and $\pi_{1_{H} \otimes 1} = \eu.$

 \end{thm}

\demo
Let $ v \in K$ and  define $\LL_{v}: A \to K$ by
\[ \LL_{v}(c) = \phi(c) v,\ \ \ c\in A. \]
Now let $a,b \in A$ and choose $a_{1}, b_{1}, \ldots , a_{n}, b_{n} \in A$
such that
\[ (*) \ \ \ \D(a) (b \otimes 1) = \sum_{i=1}^{n} b_{i} \otimes a_{i}. \]
Then observe that
\[ \sum_{i=1}^{n} \Lambda(a_i) \otimes \phi(b_i) v = (\LL \odot
\LL_{v} )(\sum_{i=1}^{n} a_{i} \otimes b_{i}) = (\LL \odot
\LL_{v} ) ( \chi (\D(a) (b \otimes 1)) ) .\]
So the left-hand side above is independent of the choice of the $a_{i}$'s
and $b_{i}$'s as long as they satisfy $(*).$ We set therefore
\[ U( \LL(a) \otimes \phi(b) v) = \sum_{i=1}^{n} \Lambda(a_i) \otimes
\phi(b_i) v .\]
Observe now that
\[ \sum_{i,j=1}^{n} b_j^*b_{i} \f(a_{j}^* a_{i}) = (\iota \otimes \f)
  ( ( \sum_{j=1}^{n} b_{j}^* \otimes a_{j}^*)( \sum_{i=1}^{n} b_{i} \otimes
  a_{i}) ) \]
  \[ = (\iota \otimes \f ) ( (\D(a) (b \otimes 1))^* (\D(a) (b \otimes
  1)) ) \]
  \[ = (\iota \otimes \f ) (  (b^* \otimes 1) \D(a^* a) (b \otimes
  1)  ) \]
  \[ = b^* b \, \f(a^* a)   \ \ \ \ \ \textnormal{(by left invariance of}
  \ \f). \]
Therefore we have
  \[ \| U ( \LL(a) \otimes \phi(b) v) \|^2 =
  \sum_{i,j=1}^{n} \ip{ \LL(a_{i}) \otimes \phi(b_{i}) v}
  {\LL(a_{j}) \otimes \phi(b_{j}) v} \]
  \[ =\sum_{i,j=1}^{n} \ip{\LL(a_{i})} {\LL(a_{j})} \ip{\phi(b_{i})
  v}{\phi(b_{j}) v } = \sum_{i,j=1}^{n} \f(a_{j}^* a_{i})
  \ip{\phi(b_{j}^*b_{i}) v}{v} \]
  \[ = \ip{ \phi( \sum_{i,j=1}^{n} b_j^*b_{i} \f(a_{j}^* a_{i})) v}
  {v} = \ip{\phi( b^* b \f(a^* a))v} {v} \]
\[  = \f(a^* a) \ip{\phi(b^*b) v}{v} =  \| \LL(a) \otimes \phi(b)
  v \|^2 \]

If we now extend $U$ by linearity to a  map from $\LL(A) \odot
\phi(A)K$ into itself, the same kind of argumentation as above shows that
$U$
is a well-defined  isometry on $\LL(A) \odot \phi(A)K,$ and thereby that the
the map  $U$  is   well
defined. Since $\phi$ is assumed to be non-degenerate, we can extend
$U$ to a linear isometry on $H \otimes K.$ Moreover, using the
cancellation properties of $\AA,$ one easily checks that $U$ is
surjective, hence that it is a unitary on $H \otimes K.$

\vspace{1ex}
 We now show that $U \in M(\har \otimes A_{\phi}).$

 \vspace{1ex}
  Let $a, b \in A.$  Write $ \D(b) (a \otimes 1) = \sum_{k=1}^{n}
a_{k} \otimes b_{k}$ for some $a_{1},b_{1}, \ldots, a_{n}, b_{n} \in A.$
Then
  \[ (**) \ \ \ \ \ U (\hat \pi (\hat b) \otimes \phi(a)) =  \sum_{k=1}^{n}
  \hat{\pi}(\hat{b}_{k}) \otimes \phi(a_{k}). \]
  To prove this, consider $c \in A, v \in K.$ It suffices to show
  that
  $$( U (\hat \pi (\hat b) \otimes \phi(a))) ( \hat{\LL}(\hat c)
  \otimes v) = (\sum_{k=1}^{n}
  \hat{\pi}(\hat{b}_{k}) \otimes \phi(a_{k})) ( \hat{\LL}(\hat c)
  \otimes v),$$
  that is,
  $$(***) \ \ \ \ \ \  U ( \hat{\LL}(\hat b \hat c) \otimes \phi(a) v) =
\sum_{k=1}^{n}
  \hat{\LL}(\hat{b_{k} \hat c}) \otimes \phi(a_{k}) v.$$
  Let L denote the left side of  equation $(***).$
  Write $b \otimes c = \sum_{i=1}^{r}
 \D(p_{i}) (q_{i}\otimes 1)$ for some $p_{1},q_{1}, \ldots, p_{r},
 q_{r} \in A.$ Then we have $\hat b \hat c =  \sum_{i=1}^{r} \varphi(q_{i})
 \hat{p}_{i}$ (as in the proof of Lemma \ref{techlem3}). Further, for
 each $i = 1, \ldots , r, $ write
 $\D (p_{i}) (a \otimes 1) =  \sum_{j=1}^{s(i)}
x_{ij} \otimes y_{ij}$ for some $x_{i1},y_{i1}, \ldots, x_{is(i)},
y_{is(i)} \in A.$

\vspace{1ex}
Then, using the definition of $U$ at the second step, we get
$$L = \sum_{i=1}^{m} \varphi (q_{i}) U ( \LL(p_{i}) \otimes \phi (a)
v )$$
$$ = \sum_{i=1}^{m} \sum_{j=1}^{s(i)}  \varphi (q_{i})  ( \LL(y_{ij})
\otimes \phi (x_{ij}) v )$$
$$= \sum_{i=1}^{m} \sum_{j=1}^{s(i)}  \hat{\LL} (\varphi (q_{i})
\hat{y}_{,j})
\otimes \phi (x_{ij}) v \, .$$
Set
$ X = \sum_{i=1}^{m} \sum_{j=1}^{s(i)}  \varphi (q_{i})   \hat{y}_{ij}
\otimes x_{ij} \, \in \hat A \odot A.$ Then, to show that $(***)$ holds, it
clearly suffices  to prove that $X=\sum_{k=1}^{n}
  \hat{b}_{k} \hat c \otimes a_{k} \, .$ But this is precisely what is
  established in Lemma \ref{techlem6}. Hence, we have shown
  that $(**)$ holds.

  \vspace{1ex}
  Let $\F:  A \to \hat A$ denote the ``Fourier" transform'', that is $\F(a)
=
  \hat a, \ a \in A. $
  Then $(**)$ may be rewritten as
  \[ U (\hat \pi (\hat b) \otimes \phi(a)) = (\hat \pi \odot \phi ) (
  \F \odot \iota) \, \chi \, (\D(b) (a \otimes 1)),\ \ \ \ a,b \in A. \]

  This means that $U ( \hat \pi (\hat A)  \odot  \phi(A) ) \subset
  \hat \pi (\hat A)  \odot  \phi(A).$ Since $\F$ is
  bijective and $\D(A) (1 \otimes A) = A \otimes A,$ we get in fact
  equality. A continuity argument  gives then
  $U ( \har \otimes A_{\phi}) = ( \har \otimes A_{\phi}).$ As $U$ is
  unitary, we also get
   $U^*( \har \otimes A_{\phi}) = $ \\ $( \har \otimes A_{\phi}).$ Applying
   the $*$-operation in $B(H \otimes K),$ we then get
    $( \har \otimes A_{\phi}) U =$ \\ $ ( \har \otimes A_{\phi}).$
    Hence, we have shown that $U \in M(\har \otimes A_{\phi}).$

\vspace{1ex} Next we show that $U \in M ( \har \otimes B_{0}(K) ).$
To this end,  we first prove that $U ( \har \otimes B_{0}(K)) \subset
\har \otimes B_{0}(K).$ Now, for all $a, b, c \in A$ and $u, v \in K,$
we have
$$ U ( \hat \pi ( \hat b \hat c ) \otimes  ( \cdot , u) \phi (a) v) =
U ( \hat \pi ( \hat b) \otimes \phi(a)) ( \hat \pi (\hat c) \otimes
(\cdot , u) v)  \ \in \hat \pi (\hat A) \odot B_{0}(K)$$
as we have seen that $U ( \hat \pi ( \hat b) \otimes \phi(a)) \in
\hat \pi (\hat A) \odot \phi(A).$ Since $\hat{A}^2 = \hat A$ and
$\phi$ is non-degenerate, it follows from a continuity argument that
 $U ( \har \otimes B_{0}(K)) \subset \har \otimes B_{0}(K),$ as desired.
 Now, using that $ U^* (\hat \pi (\hat A) \odot \phi(A)) \subset
 \hat \pi (\hat A) \odot \phi(A),$ we get similarly that
 $U^* ( \har \otimes B_{0}(K)) \subset \har \otimes B_{0}(K).$ Taking
 adjoints, we get $ ( \har \otimes B_{0}(K)) U \subset \har \otimes
 B_{0}(K).$ Hence, we have shown that $U \in M(\har \otimes B_{0}(K)).$

\vspace{1ex}

We now establish the following formula :

\[ (****) \ \ \ \ \ (\om_{\LL(a), \LL(b)} \otimes \iota ) \, U =
\phi ( ( \iota \otimes \f) (( 1 \otimes b^*) \D(a) )), \ \ \ a, b \in A. \]

Let $d \in A$ and $ v, w \in K.$ \\

\noindent Write $ \D(a) (d \otimes 1) = \sum_{i=1}^{n}
d_{i} \otimes a_{i}$ for  $a_{1},d_{1}, \ldots, a_{n}, d_{n} \in A.$
Then
$$\ip{\om_{((\LL(a),\LL(b)}\otimes \iota) \, U) \, \phi(d)  v} { w}
 = \ip{U ( \LL(a) \otimes \phi(d) v)} { \LL(b) \otimes w }$$
$$ = \sum_{i=1}^{n} \ip{\LL(a_{i}) \otimes \phi(d_{i}) v}{ \LL(b) \otimes w
}
 =  \sum_{i=1}^{n} \f(b^* a_{i}) \ip{\phi(d_{i}) v}{w}$$
$$ = \ip{\phi( \sum_{i=1}^{n} \f(b^* a_{i}) d_{i}) v} {w}
 = \ip{\phi( (\iota \otimes \f) (\sum_{i=1}^{n} d_{i} \otimes b^* a_{i})) v}
{w}$$
$$ = \ip{\phi( (\iota \otimes \f) ( (1 \otimes b^*) \D(a) (d \otimes 1)))
v} {w}$$
$$ = \ip{\phi( (\iota \otimes \f) ( (1 \otimes b^*) \D(a))) \, \phi(d) v}
{w},$$
which shows $(****).$

\vspace{2ex}
Using this formula, the norm-closure of
$ \{ (\om \otimes \iota )  \, U \, | \, \om \in B_{0}(K)^* \} $ is 
easily seen to be equal to
$A_{\phi}.$

\vspace{2ex}

We are now in position to  prove the formula relating $\phi$ and $U$, that
is,

 \[ \phi(a) = (\qr(a) \otimes \iota) \,  U , \ \ \ a \in A. \]

 Let $a \in A.$ Pick $c \in A$ such that
 $$ \hat c \widehat{S(a^*)} = \widehat{S(a^*)} \hat c =
 \widehat{S(a^*)} .$$

 Then  we have

 \[ (\qr(a) \otimes \iota) \,  U
  = (\om_{\LL(a), \LL(c)} \otimes \iota ) \, U
    \ \ \ \ \ \textnormal{(using Proposition \ref{Qmap})} \]
\[ = \phi ( ( \iota \otimes \f) (( 1 \otimes c^*) \D(a) ))
\ \ \ \ \ \textnormal{(using $(**** )$ above}) \]
\[ = \phi (a)   \ \ \ \ \ \ \ \ \ \ \ \ \ \ \ \ \textnormal{(using Lemma
\ref{techlem4})}, \]
as desired.

\vspace{2ex}
Once this fundamental formula is established, the  corepresentation property
and the uniqueness of $U$ follow readily from the norm-density of $Q(A)$
in $\hat M_{*}.$  For example, regarding $(\hdarop \otimes \iota) \, U$ and
$U_{13}U_{23}$ as lying in $\hat M \bar{\otimes} \hat M \bar{\otimes}
B(K)$, as we may, we have

\[ (Q(a) \bar{\otimes} Q(b) \bar{\otimes} \iota) ((\hdarop \bar{\otimes}
\iota) \, U)
  = ( ( (Q(a) \bar{\otimes} Q(b)) \hdarop ) \bar{\otimes} \iota) \, U \]
\[ = (Q(ab) \bar{\otimes} \iota) \, U = (\qr(ab) \otimes \iota) \, U =
\phi(ab) = \phi(a) \phi(b) \]
\[ = ((\qr(a) \otimes \iota) \, U) ((\qr(b) \otimes \iota) \, U)
  = ((Q(a) \bar{\otimes} \iota)\,  U) ((Q(b) \bar{\otimes} \iota) \, U) \]
 \[ = (Q(a) \bar{\otimes} Q(b) \bar{\otimes} \iota) (U_{13}U_{23}) \]
 for all $a,b \in A.$

 \vspace{2ex}
Finally, we  check the last assertion of the theorem.
Let $a,b,c \in A$ and choose $a_{i}$'s and  $b_{i}$'s as above.
Then we have

\[ \hat W (\LL(a) \otimes \pr (b) \LL(c)) =
   ( \Sigma W^* \Sigma) (\LL(a) \otimes  \LL(bc)) =
   ( \Sigma W^* ) (\LL(bc) \otimes  \LL(a)) \]
\[ = \Sigma (\LL \odot \LL) ( \D(a) ( bc \otimes 1)) =
   \Sigma (\LL \odot \LL) ( \sum_{i=1}^{n} b_{i}c \otimes a_{i})
   = \sum_{i=1}^{n} (\LL \odot \LL) (a_{i} \otimes b_{i} c) \]
 \[ = \sum_{i=1}^{n}  \LL(a_{i}) \otimes \pr (b_{i}) \LL(c)
    = U(\pr)( \LL(a) \otimes \pr (b) \LL(c)) \]

    and

    \[ (1_{H} \otimes 1) (\LL(a) \otimes \eu (b) 1) = \LL(a \e(b))
      = \LL ( (\e \otimes \iota ) ( \D(a) ( b \otimes 1))) \]
      \[ = \LL (\sum_{i=1}^{n} a_{i} \e(b_{i})) = \sum_{i=1}^{n}
      \LL(a_{i}) \otimes \eu(b_{i}) 1 = U(\eu) (\LL(a) \otimes \eu (b)
      1) .\]

      This  clearly implies that $\hat W = U(\pr)$ and $1_{H} \otimes
      1 = U(\eu),$ as desired. \qed

 \vspace{2ex}
\emph{Remark.}
It is clear that Theorem \ref{generator} and Theorem \ref{generator2}
together provide a bijective correspondence between unitary
corepresentations of $(\har, \hdarop)$ and
non-degenerate $*$-representations of $\au.$ In a similar way, one may
prove  that there is a bijective correspondence between unitary
corepresentations of $(\har, \hdar)$ and
non-degenerate $*$-representations of $\au,$ as proved in \cite{Ku2}
in a quite different way. Alternatively, one may use here that $U \to
U^*$ gives a bijective correspondence between unitary
corepresentations of $(\har, \hdarop)$ and  unitary
corepresentations of $(\har, \hdar).$

We also mention that Theorem \ref{generator} and Theorem \ref{generator2}
may easily be dualized to produce a bijective correspondence between unitary
corepresentations of $(\ar, \dar)$ and
non-degenerate $*$-representations of $\hat{A}_{\mathrm{u}}.$

\vspace{2ex}
\emph{Remark.}
Let $V \in \hat M \bar{\otimes} B(K)$ be a unitary such
that $(\hdarop \bar{\otimes} \iota)\, V = V_{13} V_{23}$ and let
$\tilde{\pi}_{V}$ be the associated $*$-representation of $\au$
defined in our previous remark. As a consequence of Theorem
\ref{generator2} we then get
\[(Q(a) \bar{\otimes} \iota) \, V =  \tilde{\pi}_{V}(a) \]
\[ = (\qr(a) \otimes \iota) \, U(\tilde{\pi}_{V})
= (Q(a) \bar{\otimes} \iota) \, U(\tilde{\pi}_{V}) \]
for all $a \in A.$ This implies that $V = U(\tilde{\pi}_{V}).$

Especially, we have $V \in M(\har \otimes B_{0}(K)),$ as mentioned in
a previous remark.

\vspace{2ex}
\emph{Remark.} Let $U$ be a unitary corepresentation of $(\har, \hdarop)$ on
$K.$
We can define a representation $\Phi_{U}$ of the Banach algebra
$\hat{M_{*}}$
on $K$ by
$\Phi_{U}(\om) = (\om \bar{\otimes} \iota) \,  U,$ \\
$  \om \in \hat{M_{*}}. $
Then we have $(\Phi_{U} \circ Q) (a) = \pi_{U}(a) , \ a \in A,$ so
$\Phi_{U} \circ Q$ is $*$-preserving and non-degenerate. One easily sees
that the map
$U \to \Phi_{U} $ gives a bijective correspondence between all
unitary  corepresentations of $(\har, \hdarop)$ and all
representations $\Phi$ of $\hat{M_{*}}$ such that $\Phi \circ Q$
is $*$-preserving and non-degenerate.

\section{Amenability, co-amenability and nuclearity}

We prove in this section the results stated in the Introduction.

\begin{thm} \label{univred}
    Let $\AA$ be an \aqg\ . Then $\AA$ is co-amenable if and only if
    $\au = \ar,$ that is, the canonical map $\pr$ from $\au$ onto
    $\ar$ is injective.
    \end{thm}

    \demo
    Assume that $\au = \ar.$ Then $\pr$ is an isometry.

    Let  $a \in A.$
    Then we have $\| \pu(a) \|_{\mathrm{u}} = \| \pr( \pu(a) \| = \| \pi(a) \|.$
    Hence,

    \[ | \e(a)| = | \eu(\pu(a))| \leq \| \pu(a) \|_{\mathrm{u}} =  \| \pi(a) \|. \]

    This shows that $\e$ is bounded on $A$ with respect to the reduced
    norm, that is, $(A, \D)$ is co-amenable.

    \vspace{2ex}
    Assume now that $\AA$ is co-amenable, that is,
    $| \e(a)|  \leq  \| \pi(a) \| $ for all $a \in A.$ Using
    Theorem \ref{generator2}, we can express this as
    \[ | (Q(a) \bar{\otimes} \iota) (1_{H} \otimes 1) |
    \leq | (Q(a) \bar{\otimes} \iota) \hat W | , \ \ \  a \in A. \]
    Using  the norm-density of $Q(A)$ in $\hat{M_{*}}$ (cf. Proposition
    \ref{Qmap}) and a continuity argument, we can conclude that
    \[  | \om( 1_{H}) | = |(\om \bar{\otimes} \iota) (1_{H} \otimes 1) |
    \leq \| (\om \bar{\otimes} \iota) \hat W \| , \ \ \ \om \in
    \hat{M_{*}}. \]

    To show that $\au = \ar,$
    it is enough to show that $\pr$ is isometric on $A$ = $\pu(A),$ or,
    equivalently, that $\|a\|_{\mathrm{u}} \leq \|\pi(a) \| \ , \ a \in A$
    (since the reverse inequality always holds by definition of the
    universal norm). To show this inequality,  it suffices to show that
    \[ \| \phi(a) \| \leq \| \pi(a) \| \]
    for any given non-degenerate $*$-representation $\phi$ of $A$ on some
    Hilbert space $K$ and any given $a \in A.$

    Now, let $U= U(\phi)$ be a generator for $\phi$ (extended to
    $\au$ to be pedantic), according to Theorem \ref{generator2}. Then
    this amounts to show

    \[ (*) \ \ \ \| (Q(a) \bar{\otimes} \iota ) U  \| \leq
        \| (Q(a) \bar{\otimes} \iota ) \hat W \|. \]
 
To show $(*),$ we adapt an argument from \cite[Proposition 5.5]{Bl}
(where Blanchard characterizes the amenability of regular multiplicative
unitaries). 

\vspace{2ex}
 Let $v,w \in H$, $\| v \| = \| w \| = 1.$ Define $\om \in
 \hat{M_{*}}$ by
 \[ \om ( x) = (Q(a) \bar{\otimes} \om_{v,w})((x \otimes 1) \, U), \ \ \
 x \in \hat M .\]

Then
\[ \om ( 1_{H}) = (Q(a) \bar{\otimes} \om_{v,w}) \, U = \om_{v,w}(
(Q(a) \bar{\otimes} \iota )  \, U) .\]

Hence we have

\[(**) \ \ \ | \om_{v,w}( (Q(a) \bar{\otimes} \iota ) U) | = | \om ( 1_{H})
|
\leq \| (\om \bar{\otimes} \iota) \hat W \| .\]

Now, recall (from the proof of Theorem \ref{generator}) that
we have
\[ \hat{W}_{12}^* \, U_{23} \, \hat{W}_{12} = U_{13} \, U_{23} .\]
Therefore, applying $\iota \bar{\otimes} \sigma $ to this
equation, where $\sigma$ denotes the flip map \newline from $B(H \otimes K)$ to
$B(K \otimes K),$
we get
\[ U_{32} \, \hat{W}_{13} \, U_{32}^* = \hat{W}_{13} \,  U_{12}.\]
Using this, we obtain
\[ (\om \bar{\otimes} \iota) \, \hat W =
(Q(a) \bar{\otimes} \om_{v,w} \bar{\otimes} \iota )(\hat{W}_{13}
U_{12})  \]
\[ =(Q(a) \bar{\otimes} \om_{v,w} \bar{\otimes} \iota )
(U_{32}  \hat{W}_{13}  U_{32}^* ) \]
\[ = (\om_{v,w} \bar{\otimes} \iota) ( \sigma(U)  (1_{K} \bar{\otimes} (
Q(a)
\bar{\otimes} \iota) \hat W) ) \sigma(U)^* ) ,\]
which implies that
\[(***) \ \ \ \  \| (\om \bar{\otimes} \iota) \, \hat W \| \leq \| (Q(a)
\bar{\otimes} \iota) \, \hat W \|. \]
Combining $(**)$ and $(***),$ we get
\[ | \om_{v,w}( (Q(a) \bar{\otimes} \iota ) U) | \leq
   \| (Q(a) \bar{\otimes} \iota) \hat W \| .\]
As this holds for all $v,w \in H$, $\| v \| = \| w \| = 1$, this
implies that $(*)$ holds, which finishes the proof.~\qed

\vspace{2ex}
\emph{Remark}. In \cite[Appendice]{BS}, Baaj and Skandalis
introduce the notions of amenability and co-amenability for regular
multiplicative unitaries (see also \cite{Bl}). These notions may be adapted
to
multiplicative unitaries associated to \aqg s  as follows. We first remark
that,
from the point of view adopted in \cite{BS}, it is quite natural to consider
$V =
\hat{W}$ as {\it the} multiplicative unitary associated with an \aqg\ $\AA;$
this point of view is supported by the fact that $\pi_{\hat{W}}= \pr,$
which we pointed out in Theorem \ref{generator2}. However, this is
essentially a matter of convention. The adapted Baaj-Skandalis
definition of co-amenability of  $V = \hat{W}$ amounts then to require that
$\pi_{\hat{W}} : \au \to \ar$ is injective, in which case one also says that
$W$ is
amenable. Co-amenability of $W$ and amenability of $V$ may be
defined similarily by considering $W$ to be {\it the} multiplicative unitary
associated
with $\HAA.$ Using this terminology, Theorem \ref{univred} just says
that $\AA$ is co-amenable if, and only if, $V = \hat{W}$ is
co-amenable (resp. $W$ is amenable).

\vspace{2ex}
Before stating our next result, we recall that a von Neumann algebra
$N$ acting on a Hilbert space $K$ is called injective \cite{Pa} if
there exists a linear, norm one projection map from $B(K)$ onto $N.$

\vspace{2ex}
\begin{thm} \label{injectivity}
    Assume that $\AA$ is an \aqg\ such that $(\hat A, \hat \D )$ is
    amenable. Let $\phi$ be any non-degenerate $*$-representation of
    $\au$ on some Hilbert space $K.$ Then the von Neumann algebra
    $N = \phi(\au) '' \subset B(K)$ is injective.
    \end{thm}

\demo By a  classical result of J. Tomiyama \cite{Pa,To}, we can
equivalently show that $N' =
\phi(\au)'= \phi(A)'$ is injective, that is, we have  to construct a linear,
norm one projection of
$B(K) $ onto $\phi(A)'$.

Let $U$ be a unitary generator for $\phi,$ so $U \in M(\hat \ar \otimes
B_{0}(K)) \subset \hat M \bar{\otimes} B(K),$
according to Theorem \ref{generator2}.  We introduce the
unital (injective) normal $*$-homo- \\ morphism $\alpha :  B(K) \to \hat
M \bar{\otimes} B(K) $ given by
\[ \alpha(x) = U^* (1 \otimes x) U , \ \ \ x \in B(K). \]
 Then $ \alpha$ is
 an action of $\hat M$ on $B(K),$ in the sense that $(\iota
 \bar{\otimes} \alpha) \, \alpha  = (\hdarop \bar{\otimes} \iota) \,\alpha.$
 Indeed, we have
 $$ (\iota \bar{\otimes} \alpha ) \, \alpha(x) = (\iota \bar{\otimes} \alpha
)
 (U^* (1 \otimes x)
 U) = ((\iota \bar{\otimes} \alpha )U)^* (1 \otimes \alpha(x)) ((\iota
 \bar{\otimes} \alpha
 )U)$$
 $$ =(U_{23}^* U_{13} U_{23})^* (1 \otimes U^*(1 \otimes x)U) (U_{23}^*
 U_{13} U_{23})$$
 $$= (U_{23}^* U_{13}^* U_{23}) U_{23}^* (1 \otimes 1 \otimes x)
 U_{23} (U_{23}^* U_{13} U_{23})$$
 $$= U_{23}^* U_{13}^* (1 \otimes 1 \otimes x)  U_{13} U_{23}
 = (U_{13}U_{23})^* (1 \otimes 1 \otimes x) U_{13}U_{23}$$
 $$= ((\hdarop \bar{\otimes} \iota) U)^* ( \hdarop (1) \otimes x)
 (\hdarop \bar{\otimes} \iota) U$$
 $$= (\hdarop \bar{\otimes} \iota) ( U^*
 (1 \otimes x) U) = (\hdarop \bar{\otimes} \iota) \,\alpha (x)$$
 for all $x \in B(K).$

 \vspace{1ex}
 Set
 $$B(K)^{\alpha}= \{ x \in B(K) \, | \, \alpha(x) = 1 \otimes x \}.$$
 Then, using the density of $Q(A)$ in $\hat{M}_{*}$ and Theorem
 \ref{generator2}, we have

 \vspace{2ex}
 $B(K)^{\alpha}= \{ x \in B(K) \, | \,   (1 \otimes x) U = U(1 \otimes x)
 \} $

 \vspace{1ex}
 \hspace{5ex} $= \{ x \in B(K) \, | \,
 (Q(a) \bar{\otimes} \iota)  ((1 \otimes x) U ) =(Q(a) \bar{\otimes} \iota)
(
 U(1 \otimes x)) , \, \forall a \in A \}$

 \vspace{1ex}
 \hspace{5ex} $=  \{ x \in B(K) \, | \,
 x((Q(a) \bar{\otimes} \iota)  U ) =((Q(a) \bar{\otimes} \iota) U)
 x , \, \forall a \in A \}$

  \vspace{1ex}
 \hspace{5ex} $=  \{ x \in B(K) \, | \, x \phi(a) = \phi(a) x , \, \forall a
\in A \}$

  \vspace{1ex}
 \hspace{5ex} $= \phi(A)'$.

 \vspace{2ex}
  We shall now adapt an argument of Enock and Schwartz given in the proof of
\cite[Theorem
 3.1]{ES} to  construct a linear, norm one projection from $B(K)$ onto
$B(K)^{\alpha} =
 \phi(A)'$.

 \vspace{1ex}
 Using our amenability assumption, we can pick a right-invariant
 mean $m$ for $(\hat A, \hat \D_{\mathrm{op}})$ (picking first a
 left-invariant mean for $(\hat A, \hat \D )$ and combining it with
 the anti-unitary antipode of $\hat M).$

 \vspace{1ex}
 Using that $ | m ( (\iota \bar{\otimes} \eta) \alpha(x) ) | \leq \| \eta \|
  \| x \| $ for all $\eta \in B(K)_{*} \, , \, x \in B(K),$ one easily
  sees that there
  exists  a  linear contraction map $E: B(K) \to B(K) $ such that

 \[ \eta( E(x)) = m ( (\iota \bar{\otimes} \eta) \alpha(x) ) \, , \ \
 \ \eta \in
 B(K)_{*} \, , \, x \in B(K). \]

 For $\om \in \hat{M}_{*}  \, , \, \eta \in B(K)_{*}\, ,$ we have
 $$ (\om \bar{\otimes} \eta) ( \alpha (E(x)) = ((\om \bar{\otimes} \eta)
\circ \alpha)
 (E(x))$$
 $$=  m ( (\iota \bar{\otimes} (\om \bar{\otimes} \eta) \circ \alpha)
\alpha(x) )
 = m ( (\iota \bar{\otimes} \om \bar{\otimes} \eta) ( \iota \bar{\otimes}
\alpha)
 \alpha(x) )$$
 $$=  m ( (\iota \bar{\otimes} \om \bar{\otimes} \eta) (\hdarop
\bar{\otimes} \iota) \,
 \alpha (x))
 = m ( (\iota \bar{\otimes} \om) \hdarop ( ( \iota \bar{\otimes} \eta)
 \alpha(x)) )$$
 $$= m ( ( \iota \bar{\otimes} \eta) \alpha(x) )
 \om(1) \ \ \ \ \ \ \ \ \ \ ( \textnormal{ using
 right-invariance of} \ m) $$
 $$= \eta (E(x)) \om(1) = (\om \bar{\otimes} \eta) (1 \bar{\otimes} E(x)) \, .$$
 It follows that $ \alpha (E(x)) = 1 \otimes E(x) $ for all $ x \in
 B(K),$ hence that $E$ maps $B(K)$ into $B(K)^{\alpha}.$
 Further, if $x \in B(K)^{\alpha},$ that is, $\alpha(x) = 1 \otimes
 x,$ then
 $$ \eta(E(x)) = m ( (\iota \bar{\otimes} \eta) (1 \otimes x) ) =
 m(1) \eta(x) = \eta(x)$$
 for all $\eta \in B(K)_{*}.$ Thus, $E(x) = x$ for all $x \in
 B(K)^{\alpha}.$ It clearly follows that $E$ is a norm one projection
 from $B(K)$ onto $B(K)^{\alpha},$ which finishes the proof.~\qed

\vspace{2ex}
\begin{cor} \label{nuclearity} Assume that $\AA$ is an \aqg\ such that
$(\hat A, \hat \D )$ is
    amenable. Then $\au$ is nuclear.
\end{cor}
\demo By applying Theorem \ref{injectivity} to the universal
$*$-representation $\Phi$ of $\au,$ we obtain that  the second dual
$\au^{**} =
\Phi(\au)''$ is injective. By a famous result of Connes, Choi and
Effros (see \cite[2.35]{Pa} for references), this is  equivalent to
the nuclearity of $\au.$~\qed

\vspace{2ex} We shall now give  a simplified proof of a result which is
essentially due to
Ruan (see \cite[Theorem 4.5]{Ruan}).

\begin{thm} \label{Ruan}
Assume that $\AA$ is an \aqg\ such that its associated von Neumann algebra
$M \subset B(H) $ is
injective. Assume further that $\AA$ is compact with unit 1 and has a
tracial Haar
functional (that is, equivalently, $(M, \dar)$ is a compact Kac
algebra \cite{ES}).

\vspace{1ex}
Then $\AA$ is co-amenable.
\end{thm}

\demo  As usual in the compact case, we  work with the
normalized Haar functional $\f$ of $\AA.$ It is known \cite{BS,Wo} that the
traciality of
$\f$ is equivalent to
 $S^2 = \iota,$ or, equivalently, to $S$ being $*$-preserving.

\vspace{1ex}
Using the traciality assumption and the fact that $\f$ is
$S$-invariant, it  is  straightforward to check that the linear map $V_{0}:
\LL(A) \to \LL(A)$
defined by $V_{0}(\LL(a)) = \LL(S(a)), \, a \in A,$ is an isometry,
which extends to a self-adjoint unitary $V$ on $H.$ (See  \cite[Proposition
5.2]{BS} for a similar statement in the non-tracial case, which we
will use in the proof of our next result).  A simple calculation gives
$V \pi(a) V \pi(b) = \pi(b) V \pi(a) V $ for all $a, b \in A.$
Hence, $\Ad(V)$ maps $\ar$ (and  $M$) into $ \pi(A)' = M'.$

\vspace{1ex}

Now, we recall that  injectivity of $M$
implies that the $*$-homomorphism \\ $P : M \odot M' \to B(H)$
determined by $$P(x \otimes y) = x \, y , \, x \in M, \, y \in M',$$
has a
bounded extension  $\tilde{P}: M \otimes M' \to B(H), $
where we stress
that $\otimes$ denotes the minimal tensor product (as opposed to the von
Neumann
algebra tensor product). (Note that this deep result is not mentioned
explicitely in
\cite{Pa}. It may be deduced from the literature as follows.
Injectivity is
equivalent to semidiscreteness, as first shown by Connes \cite{Co} in the
factor
case. A  direct  proof of the forward
implication due to Wassermann may be found in  \cite{W}. The backward
implication is shown by Effros and Lance \cite{EL},
who also show that semidiscreteness is equivalent to the above
property).

\vspace{1ex}
We use $\tilde{P}  $ to define a map $\e_{0} : \ar \to \bC$ by
\[ \e_{0}(x) = \ip { ((\tilde{P} \circ  ( \iota \otimes \Ad (V) ) \circ
\dar) (x) ) \, \LL(1)} { \LL(1)}, \ \ x \in \ar. \]

Clearly, $\e_{0}$ is a  state on $\ar.$ Further,
we have $\e_{0}(\pi(a)) = \e(a), a \in A.$ Before establishing this fact, we
point out  that it clearly implies that $\e$ is bounded with respect to
the reduced norm of $\ar,$ that is,  $\AA$
is co-amenable.

\vspace{1ex}
Let $a \in A$ and write $ \D(a) = \sum_{i=1}^{n} a_{i} \otimes b_{i}$
for some $a_{1}, b_{1}, \ldots ,a_{n}, b_{n} \in A.$ Then
$$((\tilde{P} \circ  ( \iota \otimes \Ad (V) ) \circ
\dar) (\pi(a)) ) \, \LL(1)$$
$$ = ((\tilde{P} \circ  ( \iota \otimes \Ad (V) ) \circ (\pi \odot
\pi) (\D(a)) ) \, \LL(1)$$
$$ = (\sum_{i=1}^{n} \pi(a_{i}) \Ad(V)(\pi(b_{i})) \, \LL(1)$$
$$ = \LL ( \sum_{i=1}^{n} a_{i} S(b_{i}) ) = \LL (m (\iota \odot S)
\D(a) ) $$
$$ = \LL( \e(a) 1) = \e(a) \LL(1).$$
Hence,
$$ \e_{0}(\pi(a)) = \ip { ((\tilde{P} \circ  ( \iota \otimes \Ad (V) ) \circ
\dar) (\pi(a))  \, \LL(1)} { \LL(1)}$$
$$ = \ip {\e(a) \LL(1)} {\LL(1)} = \e(a),$$
as asserted.~\qed

\vspace{2ex}
It would be interesting to know whether the traciality assumption
in Theorem \ref{Ruan} is redundant. We shall now prove a related result,
saying that injectivity of $M$ implies  a  kind of perturbed
co-amenability of $\AA.$  We  recall first some more facts about the
compact case.

\vspace{1ex}
Let  $\AA$ be a compact \aqg\ with unit $1.$ It is immediate that $(\ar,
\dar)$ is a
compact quantum group in the sense of Woronowicz \cite{Wo1,Wo},
with Haar state $h_{r}$  given by the restriction of $\om_{\LL(1)}$ to
$\ar.$
The  unique dense Hopf $*$-subalgebra \cite{BMT} of $(\ar, \dar)$ may
 be identified with
 $(A, \D, \e, S)$ ( via the Hopf $*$-algebra isomorphism $\pi$).
 Using this identification, we may introduce the remarkable family
 $(f_{z})_{z \in \bC}$ of multiplicative linear functionals on $A$
 constructed by Woronowicz (see \cite{Wo1, Wo}).

 \vspace{1ex}
 Some of the properties of this family are  $f_{0}= \e;$ $f_{z} *
 f_{z'} = f_{z+z'},$ where $\om * \eta = (\om \otimes \eta) \D, \ \
 \om , \eta \in A';$ the maps $a \to f_{z} * a = (\iota \otimes f_{z})
 \D (a) $ and $ a \to (f_{z} \otimes \iota) \D (a) $ are automorphisms
 of $A;$ we have $f_{z}^* = f_{-\bar{z}}$ and $ f_{z} \circ S =
 f_{-z};$ for all $a, b \in A,$ we have $\f (a b) = \f ( b
 (f_{1}*a*f_{1}))$ and $S^2 (a) = f_{-1}* a * f_{1}.$

 \vspace{1ex}
 \noindent We also  mention that the following three conditions are
equivalent:

 \vspace{1ex}
 $\f$ is tracial; $f_{z} = \e $ for all $z \in \bC;$ $f_{1}= \e.$

  \vspace{1ex}

\begin{thm} \label{qdim}
 Assume that $\AA$ is a compact \aqg\ such that its associated von Neumann
algebra
 $M$ is injective. Let $(u^{\alpha})$ denote a complete set of
 pairwise inequivalent irreducible  unitary
 corepresentations of the compact quantum group $(\ar, \dar)$ and let
 $n_{\alpha}$ (resp. $d_{\alpha}$)  denote the ordinary (resp.  quantum)
 dimension of $u^{\alpha}.$

 Then there exists a state $\e_{1}$ on
 $\ar$ such that
 \[ \e_{1}(u_{ij}^{\alpha}) =
 \frac{n_{\alpha}}{d_{\alpha}} \, \d_{ij} , \ \ \  1
 \leq i, j \leq n_{\alpha}  \]
 \end{thm}

 \demo We recall first that $d_{\alpha}= \sum_{i=1}^{n_{\alpha}}
f_{1}(u_{ii}^{\alpha}).$
 In other words,  $d_{\alpha}$ is the trace of the matrix
 $F_{\alpha}= (f_{1} \odot \iota ) u^{\alpha}.$

 \vspace{1ex}
 Next, we define  $\tilde{S}: A \to A$ by $\tilde{S}(a) = f_{1} *
 S(a), \ a \in A.$ (This map is sometimes called the twisted antipode of
$\AA$).
 Using the properties of the $f_{z}$'s mentioned above, one checks easily
 that $\tilde{S}$ is an involutive anti-automorphism of $A.$

 \vspace{1ex}
 Further, as shown in the proof of \cite[Proposition 5.2]{BS},  the linear
map $U_{0}: \LL(A) \to \LL(A)$
defined by $U_{0}(\LL(a)) = \LL(\tilde{S}(a)), a \in A,$ is an isometry,
which
extends to a self-adjoint unitary $U$ on $H$  satisfying
$$ (\Ad(U) (\pi(a))) \LL(b) = \LL(b \tilde{S}(a)), \ \ \ a, b \in A.$$

It follows readily that $\Ad (U) (x) \in  \pi(A)' = M',$ for all $ x \in
M.$ (In fact, one may check that $\Ad (U) (x) = J R(x^*) J$  for all $ x \in
M,$ where $J : H \to H$ denotes the Tomita-Takesaki map such that $(H,J)$ is
standard for $M.$)

\vspace{1ex}
Now, let
$\tilde{P}$ be as in the  proof of Theorem \ref{Ruan}. We then
define a state $\e_{1}$ on $\ar$ by
\[ \e_{1}(x) = \ip { ((\tilde{P} \circ ( \iota \otimes \Ad (U) ) \circ
\dar) (x) ) \, \LL(1)} {
\LL(1)}, \ \ \  x \in \ar. \]
Using the orthogonality relations \cite{Wo1, Wo} for the
$u_{ij}^{\alpha}$'s, one checks  that $\e_{1}$ satisfies the stated
property. More precisely, the computation goes as follows.
Fix $\alpha,$ set $n=n_{\alpha}, d=d_{\alpha}$ and write $u_{ij}^{\alpha} =
\pi
(v_{ij}), \, v_{ij} \in A.$ The matrix $(v_{ij})$ is then an
$n \times n$ unitary matrix over $A,$ and one of the orthogonality
relations for the $u_{ij}^{\alpha}$'s gives
$$ \f(v_{ik}v_{js}^{*}) = (1/d) \, \d_{ij}
f_{1}(v_{sk}), \ \ \ i, j, k, s \in \{1, \ldots , n\}.$$
Using this, we get

$$ \e_{1}(u_{ij}^{\alpha}) = \ip { ((\tilde{P} \circ ( \iota \otimes \Ad (U)
) \circ
(\pi \odot \pi) (\D(v_{ij}) ) \, \LL(1)} { \LL(1)}$$
$$= \sum_{k=1}^{n} \ip { \pi(v_{ik}) (\Ad(U)(\pi(v_{kj}))) \, \LL(1)}
{\LL(1)}$$
$$= \sum_{k=1}^{n} \ip { \LL(v_{ik} \tilde{S}(v_{kj}) )} {\LL(1)}
  = \sum_{k=1}^{n} \f (v_{ik} \tilde{S}(v_{kj}))$$
$$=\sum_{k=1}^{n} \f (v_{ik} (f_{1}* v_{jk}^{*}))
  =\sum_{k=1}^{n} \f (v_{ik} ((\iota \odot f_{1}) \D(v_{jk}^{*}) ))$$
$$=\sum_{k,s=1}^{n} \f(v_{ik}v_{js}^{*}) f_{1}(v_{sk}^{*})
  = (1/d) \sum_{k,s=1}^{n} \delta_{ij} f_{1}(v_{sk}) f_{1}(v_{sk}^{*})$$
$$= (1/d) \, \delta_{ij} \, f_{1} (\sum_{k,s=1}^{n} v_{sk}v_{sk}^{*})
  = (1/d) \, \delta_{ij} \, f_{1} ( n\, 1)  = (n/d) \, \delta_{ij}$$
as desired.~\qed

\vspace{2ex}
\emph{Remark.} Assume the existence of a state $\e_{1}$ satisfying the
statement  of
Theorem \ref{qdim}. If we
 also assume that $d_{\alpha} = n_{\alpha}$ for all $\alpha$
 (especially, if we also assume that
 $\AA$ has a tracial Haar functional),
 then $\e_{1} \circ \pi$  coincides with the co-unit $\e$ of $\AA$, and we
 can then conclude that $\AA$ is co-amenable. Hence, Theorem
 \ref{Ruan} is  just a special case of Theorem \ref{qdim}.

In the general case, it is known that the ordinary dimension is always
smaller than
the quantum dimension (which is always positive and less than infinity).
Thus $q_{\alpha} = \frac{n_{\alpha}}{d_{\alpha}} \in <0, 1]$. However, for
the relevant examples (like
quantum $SU(2)$ etc) it tends exponentially to zero with 'increasing'
$\alpha$'s.
Of course, one may wonder whether it is possible to use the existence of
the state $\e_{1}$ to deduce that  $\e$ is bounded.

One  natural way to proceed is  to consider $\e_{1}$ as an element of
the Banach algebra $A_r^*$ and try to use spectral calculus to ``press'' up
the values
$q^{\alpha}$ to $1$. For any function $f$ which is analytic on a
region in the complex plane containing the closed unit disk with
center at the origin and satisfies $f(0) = 0$, one may check that
$$f(\e_{1} ) u^{\alpha}_{ij}=f(q_{\alpha})\delta_{ij},$$
for all $\alpha$ and $i,j$. However, it seems difficult to  proceed
further without introducing some other assumptions.
We also mention that  $A_r^*$ is  in fact a Banach $*$-algebra with
$*$-operation given by
$$\psi ^*(a)=\overline{\psi (R (a^* ))},$$
for all $\psi\in A_r^*$ and $a\in A$.
One may  show that $\e_{1} R = \e_{1}$, hence that $\e_{1}^{*}=
\e_{1},$ but it is not clear that this may be of any help.

Another possible approach is to consider the  bounded linear map from
$A_r$ into itself given by
$\psi_{1} =(\iota \otimes \e_{1} )\D_r \, .$
It is not difficult to show that it is injective. If one could show
surjectivity
of $\psi_{1},$ then, appealing to the Open Mapping Theorem,
$\psi_{1}^{-1}$ would be bounded. Further, we would have $\e \circ
\pi^{-1}  = \e_{1}
\circ \psi_{1}^{-1}$ on $\pi (A).$ Hence, we would then be able to
conclude that $\AA$ is co-amenable.
We are so far  only able to see that $\psi_{1}$ has dense range, as it
contains $\pi (A).$

\vspace{2ex}
\emph{Remark.} Let $J$ and $\nabla $ be the Tomita-Takesaki maps
associated to $\LL(A),$ considered as  a left Hilbert algebra in $H,$
so that
$$J\nabla^{1/2}\Lambda (a)=\Lambda (a^* ),$$
for all $a\in A$. Further, let $\tau$ denote the scaling group of
$\AA$ (see \cite{KD}).

For $z\in {\bf C},$ define a map $Q_z :\pi
(A)\rightarrow M'$ by
$$Q_z \pi (a)=J\nabla^{-iz}(R\tau_{-z}\pi
(a))^*\nabla^{iz}J,$$
for all $a\in A$.

Clearly, $Q_z$ is  unital, muliplicative and linear. Setting $z =
i/2$ gives,
$$Q_{i/2}\pi (a)=J\D^{1/2}\pi (S(a)^*)\D^{-1/2}J,$$
for all $a\in A$.

If $z=t$ is real, then $Q_{t}$ is $*$-preserving and bounded, and it may
be extended to $M$. Note also that
$$Q_0 (x)=JR(x^* )J, \ \  x \in M.$$
Now, define a unital linear functional $\phi_z$ on $\pi
(A)$ by
$$\phi_z \pi (a)=(P(\pi \odot Q_z \pi )\D (a)\Lambda
(1), \Lambda (1)),$$
for all $a\in A$ and $z\in {\bf C}, P : M \odot M' \to B(H)$ being defined
as in the proof
of Theorem \ref{Ruan}.

For general $z\in {\bf C}$, one may easily show that
$$\phi_z \pi (u^{\alpha}_{lj})
=\sum_k \varphi (u^{\alpha}_{lk} f_{1+i(z-\bar{z})}*
((u^{\alpha}_{jk})^* )*f_{-i(z+\bar{z})}),$$
for all $\alpha$ and all $l,j$.

From this we see that
$$\phi_{i/2}\pi (u^{\alpha}_{ij})=\delta_{ij}=\varepsilon
(u^{\alpha}_{ij}),$$
for all $\alpha$ and $i,j$, so $\phi_{i/2}\pi =\varepsilon.$ Hence,
co-amenability of $\AA$ is equivalent to the boundedness of
$\phi_{i/2}.$

Now, observe that when $z=t$ is real and $M$ is assumed to be
injective, then $\phi_t$ may be extended to a state on $\ar$ such that
$$\phi_t (x)=(\tilde{P}(\iota \otimes Q_t )\D_r (x)\Lambda (1),
\Lambda (1)),$$
for all $x\in A_r, \ \tilde{P}$ being defined as in the proof of
Theorem \ref{qdim}. Note that $\phi_{0}$ is then just equal to the state
$\e_{1}$ obtained in this theorem. One may wonder whether some
analytic continuation argument could be used in this situation to deduce
that
$\phi_{i/2}$ is
bounded.

Co-amenability of $\AA$ may be characterized by the existence of a
non-zero
multiplicative linear functional on $\ar$ \cite{BMT2}. However, when
$\ar=SU_{q}(2), q \in (0,1),$ we have checked that none of the $\phi_{t}$
are multiplicative, even though $SU_{q}(2)$ is known to be co-amenable.

 \vspace{2ex}
\emph{Remark.} Some of the essence of Theorems \ref{Ruan} and
\ref{qdim} may be presented in a  more conceptual manner.
Assume that $\AA$ is a compact \aqg. We define the {\em adjoint}
representation $C$ of $A$ on $B(H)$ as follows. Let $P$ be the map
introduced in the proof of Theorem \ref{Ruan} and $U$ be the unitary
on $H$ introduced in the proof of Theorem \ref{qdim}. Then set
$$ C(a) = (( P \circ (\iota \odot \Ad(U)) \circ \dar)(\pi(a)), \ \ a \in
A.$$
(A more explicit way of defining  $C$ is
$$ C(a) \LL(b) = \sum_{i=1}^{n} \LL ( a_{i} b \tilde{S} (a_{i}')),$$
for $a, b \in A,$ and $\D(a) = \sum_{i=1}^{n} a_{i} \otimes a_{i}'$).

\vspace{1ex}
Using the map $\tilde{P}$ introduced in the proof of Theorem
\ref{Ruan},  one easily deduces that the injectivity of $M$
implies that $C$ is weakly contained in
$\pi,$ that is, more precisely, that  the associated
$*$-representation $C_{u}$ of $\au$ is weakly
contained in $\pr.$ On the other
hand, if the Haar state of $\AA$ is tracial, then $\tilde{S}=S$ and
the last part of
the proof of Theorem \ref{Ruan} shows that
 $\e_{u}$ is weakly contained in $C_{u}.$  Combining these two
assertions reproves Theorem \ref{Ruan}. An open question
is then whether $\e_{u}$ is always weakly contained in $C_{u}.$
A negative answer to this question is not unlikely, and
it would  then be of interest to find a  more general condition than
traciality of the Haar state  ensuring the weak containment of
$\e_{u}$ in $C_{u}$.
 
\vspace{2ex}
\noindent We conclude with a proof of Theorem \ref{main}.

\vspace{1ex}
\noindent \emph{Proof of Theorem \ref{main}}:

\vspace{1ex}
(1) implies (2): This result is shown in \cite{BMT2}.

\vspace{1ex}
(2) implies (3):  This is Corollary \ref{nuclearity}.

\vspace{1ex}
(3) implies (4): As $\ar$ is a quotient of $\au,$ this follows from
the fact that  quotients of  nuclear C$^*$-algebras are nuclear
\cite{EL}.

\vspace{1ex}
(4) implies (5):  As $M = \ar '',$ this follows from the fact that any
von Neumann algebra generated by a nuclear C$^*$-algebra  is injective
(this is easily seen by using that the double dual of a
nuclear C$^*$-algebra is injective,
as pointed out  in the proof of Theorem \ref{injectivity}).

\vspace{1ex}
Assume that $\AA$ is compact and has a tracial Haar functional.
Then (5) implies (1) is shown in Theorem \ref{Ruan}.~\qed

\vspace{2ex}
Finally, we remark that different proofs of (1) implies (5), and
of (5) implies (2) in the compact tracial case, were given in \cite{BMT2}.

\bigskip
{\parindent=0pt \parskip=4pt Addresses of the authors:

\smallskip Erik B\'edos, Institute of Mathematics, University of
Oslo, \\ P.B. 1053 Blindern, 0316 Oslo, Norway. E-mail:
bedos@math.uio.no

Gerard J. Murphy, Department of Mathematics, National
University of Ireland, Cork, Ireland. E-mail: gjm@ucc.ie

Lars Tuset, Faculty of Engineering, Oslo University College,
\\ Cort Adelers Gate~30, 0254 Oslo, Norway.
E-mail: Lars.Tuset@iu.hio.no}

\end{document}